\documentclass[11pt]{amsart}
\usepackage[T1]{fontenc}
\usepackage{lmodern,amsmath,amssymb,amsthm,mathtools,geometry,microtype,enumitem,booktabs,needspace,hyperref}
\hypersetup{colorlinks=true,linkcolor=blue,citecolor=blue,urlcolor=blue,
 pdfauthor={Evgeny Abakumov, Isabelle Chalendar, Alexei Poltoratski},pdfsubject={Zero sets, exactness, Clark measures, and analytic coefficient spaces}}
\hypersetup{bookmarksdepth=2}
\setlist[enumerate]{leftmargin=*,itemsep=3pt,topsep=5pt}
\setlist[itemize]{leftmargin=*,itemsep=3pt,topsep=5pt}
\newtheorem{theorem}{Theorem}[section]
\newtheorem{maintheorem}{Theorem}

\newtheorem{lemma}[theorem]{Lemma}
\newtheorem{proposition}[theorem]{Proposition}
\newtheorem{corollary}[theorem]{Corollary}
\theoremstyle{definition}
\theoremstyle{remark}\newtheorem{remark}[theorem]{Remark}
\numberwithin{equation}{section}
\newcommand{\D}{\mathbb D}
\newcommand{\T}{\mathbb T}
\newcommand{\lp}{\ell_A^p}
\newcommand{\E}{\mathbb E}

\newcommand{\mT}{m_{\T}}
\newcommand{\dm}{\,d\mT}

\title[Zero sets with $p$-summable Taylor coefficients]{Zero sets in spaces of analytic functions with $p$-summable Taylor coefficients}
\hypersetup{pdftitle={Zero sets in spaces of analytic functions with p-summable Taylor coefficients}}

\date{}
\begin{document}

		\author{Evgeny Abakumov}
	\address{Université Gustave Eiffel, LAMA, (UMR 8050), UPEM, UPEC, CNRS, 77454, Marne-la-Vallée, France}
	\email{evgueni.abakoumov@univ-eiffel.fr}

	\author{Isabelle Chalendar}
	\address{Université Gustave Eiffel, LAMA, (UMR 8050), UPEM, UPEC, CNRS, 77454, Marne-la-Vallée, France}
	\email{isabelle.chalendar@univ-eiffel.fr}
	\thanks{This research is partly funded by the Fédération de Recherche d’Animation Bézout.}
	
\author{Alexei Poltoratski}
\address{University of Wisconsin\\ Department of Mathematics\\ Van Vleck Hall\\
	480 Lincoln Drive\\
	Madison, WI  53706\\ USA }
\email{poltoratski@wisc.edu}
\thanks{The third author was partially supported by NSF Grant DMS-2244801 and also received support from the Simons Fellows in Mathematics program for 2026--27.}

\begin{abstract}
A zero set for a space of analytic functions is a set on which some nonzero
function in the space vanishes; it is exact if it coincides with the full zero
set of such a function. We construct zero sets that are not exact in the
classical coefficient spaces $\ell_A^p$ of functions analytic on the unit disk
for every $1\le p<2$. For $1<p<2$, a single Blaschke set $W$ is a zero set for
every $\ell_A^p$ in this range, but is not exact for any of them. The endpoint
case $p=1$, corresponding to the analytic Wiener algebra $\ell_A^1$, requires
a separate construction, yielding a non-exact zero set that accumulates at
only one point of the unit circle. The proofs combine quadratic characters,
finite-dimensional Clark transforms, and logarithmic estimates in zero-free
regions. We also establish sharp moment duality and an exact summability
criterion for a family of regular polygons. For $p>2$, we construct two
disjoint exact zero sets whose union is a uniqueness set. For these exponents,
suitably chosen rotations also yield a continuum-sized family of nonzero
closed shift-invariant subspaces of $\ell^p$ with pairwise trivial
intersections, recovering a known existence result.
\end{abstract}
\maketitle
\tableofcontents

\section{Introduction}\label{sec:intro}
\subsection{From coefficient norms to zero sets}
Let $\D=\{z\in\mathbb C:|z|<1\}$ be the open unit disk, and let
$\T=\partial\D$ be the unit circle. We work with the analytic coefficient
spaces on the unit disk
\[
 \ell_A^p=\left\{f(z)=\sum_{n\ge0}a_nz^n:
       \|f\|_p:=\left(\sum_{n\ge0}|a_n|^p\right)^{1/p}<\infty\right\},
       \qquad 1\le p<\infty.
\]
The space $\ell_A^2$ is the classical Hardy space $H^2$, and
$\ell_A^1$ is the analytic Wiener algebra. For conjugate exponents we write
$q=p/(p-1)$ when $1<p<\infty$, and set $q=\infty$ when $p=1$.

In his celebrated 1949 paper on two problems concerning linear transformations in
Hilbert space~\cite{Beu}, Beurling characterized the closed subspaces of the
Hardy space $H^2$ (identified with $\ell_A^2$ via Taylor coefficients) that are
invariant under multiplication by $z$. Since multiplication by $z$ corresponds
to the right shift on Taylor coefficients, we call such subspaces
\emph{shift-invariant}. His theorem illustrates a fruitful meeting point of
function theory and operator theory: the structure of shift-invariant
subspaces is described through the analytic properties of functions on the
disk. As his concluding remarks indicate, Beurling was already interested in
how this picture might change for spaces defined by other coefficient norms,
including $\ell_A^p$ when $p\ne2$.

The situation is markedly different for $p>2$: $\ell_A^p$ has nonzero
closed shift-invariant subspaces with pairwise trivial intersections~\cite{Ab},
and shift-invariant subspaces of every index
$n\in\{1,2,\ldots\}\cup\{\infty\}$~\cite{AB}. Both phenomena are absent
for $p=2$: any two nonzero shift-invariant subspaces of $H^2$ have a
nontrivial intersection, and every nonzero shift-invariant subspace of $H^2$
has index one, by Beurling's theorem.

Related structural phenomena have been studied in $\ell_A^p$ for both
$1<p<2$ and $p>2$. Although many interesting results are known, the theory of
the coefficient spaces $\ell_A^p$ remains at a comparatively early stage of
development, especially in contrast with the extensive theories of other
classical spaces of analytic functions on the unit disk, such as the Bergman
and Dirichlet spaces. A comprehensive account of much of the theory known at the time
is given in~\cite{CMRbook}.

Beurling's work also highlights a broader methodological point: the theory
of analytic functions, and in particular the study of their zeros, can serve
as a tool in problems about cyclic vectors and invariant subspaces of
operators. This provides a natural context for examining zero sets in
coefficient spaces, without identifying the zero-set problem with the
invariant-subspace problem itself.

For a nonzero analytic function $f$ in $\D$, write
\[
 Z(f):=\{z\in\D:f(z)=0\},
\]
where zeros are regarded as points, without multiplicity. Thus exactness below
is a property of the underlying set; multiplicities will be specified
separately when a result concerns them.
Given a space $X$ of analytic functions, we call $W\subseteq\D$ a
\emph{zero set} for $X$ if $W\subseteq Z(f)$ for some nonzero
$f\in X$, and an \emph{exact zero set} for $X$ if $W=Z(g)$ for some
nonzero $g\in X$. The distinction is important: the first condition asks
for a function vanishing on every point of $W$, while the second requires
that it have no additional zeros. Since the zeros of a nonzero analytic
function are isolated, every zero set is discrete in $\D$, in the sense
that its intersection with each compact subset of $\D$ is finite.

For $H^2$, the zero-set problem has a particularly clean answer. A discrete set $W\subset\D$ is the zero set of a
nonzero function in $H^2$ if and only if it satisfies the
\emph{Blaschke condition}
\[
 \sum_{w\in W}(1-|w|)<\infty.
\]
The associated Blaschke product realizes such a set as an exact zero set.
Thus, in $H^2$, the notions of zero set and exact zero set coincide. This characterization belongs to the classical
Hardy-space theory and is closely connected with inner--outer factorization
and Beurling's description of the closed shift-invariant subspaces of $H^2$;
see, for example,~\cite{Du,Nik}.

Compared with the classical theories of the Hardy, Bergman, and Dirichlet
spaces, relatively little is known about zero sets in $\ell_A^p$ for
$p\ne2$, and the available results remain rather fragmentary. One early
result is that, for $p>2$, there exist zero sets for $\ell_A^p$ that do
not satisfy the Blaschke condition~\cite{Vin}.

For $p\ge2$, every Blaschke sequence is the exact zero set of a function
in $\ell_A^p$: its Blaschke product belongs to
$H^\infty\subset H^2\subset\ell_A^p$, where $H^\infty$ denotes the
bounded analytic functions on $\D$. For $1\le p<2$, the inclusion goes
in the opposite direction, $\ell_A^p\subset H^2$, so every zero set
must satisfy the Blaschke condition. Nevertheless, that condition is no
longer sufficient. A theory of generalized inner functions was developed,
and vanishing sets were characterized in terms of boundedness of finite
extremal functions, namely functions minimizing the coefficient norm under
finitely many prescribed zero conditions and a normalization at zero~\cite{CMR}.
The failure of sufficiency of the Blaschke condition for $1<p<2$ was established
in~\cite{CD}, and nesting of the zero-set classes as the exponent varies was
established in~\cite{C}. A related, though distinct,
structural question is the failure of canonical factorization in $\ell_A^p$~\cite{CDfactor}. For further results and references, we refer the reader to~\cite{CMRbook}.

\subsection{Zero sets and exact realization}
The question is whether a zero set can fail to be exact: can every function
vanishing on $W$ be forced to have additional zeros? A function
$f\in X\setminus\{0\}$ with $W\subset Z(f)$ is called a
\emph{vanishing witness}, or simply a \emph{witness}, for $W$; if
$W=Z(g)$, then $g$ is an \emph{exact witness}. A witness is
\emph{normalized} when its value at zero is one; this normalization is
used only for sets not containing zero. A \emph{uniqueness set} is a set
on which no nonzero member of $X$ vanishes.

The preceding Blaschke-product argument settles exact realization in
$H^2$, but does not answer the question for $\ell_A^p$ when $1\le p<2$.
Nor does it settle exactness of non-Blaschke zero sets when $p>2$.
The following theorem answers this question for every $1\le p<2$.
\begin{maintheorem}[Non-exact zero sets in the disk]\label{main:nonexact}
\begin{enumerate}
\item There is a Blaschke set $W$ that is a zero set for every $\ell_A^p$, $p>1$,
but is not an exact zero set for any $1<p<2$. More precisely, there is a
function $F\in\bigcap_{p>1}\ell_A^p$, with $F(0)=1$ and only simple zeros,
such that $W\subset Z(F)$. The same set $W$ is exact for every $p\ge2$.
\item There exist $f\in\ell_A^1\setminus\{0\}$ and a set
$W^{(1)}\subset Z(f)$ whose only boundary accumulation point is $1$,
such that $W^{(1)}$ is not an exact zero set for $\ell_A^1$.
Nevertheless $W^{(1)}$
is an exact zero set for $H^\infty$ and for the disk algebra
$A(\D)$, the analytic functions on $\D$ continuous on $\overline\D$.
\end{enumerate}
\end{maintheorem}

The two assertions are proved in Theorems~\ref{thm:localnonexact}
and~\ref{thm:endpoint}, respectively.

In $H^2$, if $W$ is a Blaschke set and $B_W$ is its Blaschke product,
then the subspace of functions vanishing on $W$ is $B_WH^2$; in
particular, $B_W$ is an exact witness. In contrast, Theorem~\ref{main:nonexact} gives, for
each $1\le p<2$, a set $W$ such that every nonzero function in $\ell_A^p$
vanishing on $W$ has at least one additional zero. These additional zeros
need not be common to all such functions; the point is that no exact
witness exists.

The claim concerns zeros in the open disk. To our knowledge, these are the
first examples of non-exact zero sets in classical disk spaces. In the
familiar disk spaces the two notions coincide: Blaschke products give exact
realization in Hardy spaces~\cite{Du,Nik}; contractive zero-divisors give
it in the standard Bergman spaces $A^r$, $1\le r<\infty$~\cite{DKSS};
and a representation formula for the Dirichlet integral shows that
removing Blaschke factors preserves the classical Dirichlet space~\cite{Car}.
Here a contractive zero-divisor is a function that realizes prescribed
zeros and can be divided out without increasing the norm.
Indeed, if a Dirichlet function vanishes on $W$, one may divide it by the
Blaschke product formed from all its zeros outside $W$, with their
multiplicities, to obtain a Dirichlet function whose zero set is exactly
$W$. The zero-set characterizations for analytic functions smooth on the
closed disk and for analytic Lipschitz classes of positive exponent likewise
imply exact realization of subsets of their zero sets~\cite{TW,TWzeros}.
These comparisons concern the standard spaces just listed.

Such phenomena are well known in spaces of entire functions, such as
Paley--Wiener and Fock spaces. They can already be seen from classical necessary conditions on the
zero distributions of entire functions of finite type; see, for example,
~\cite[Theorem~2.10.1]{Boas}.\footnote{For example, $\mathbb N$ is contained in the zero set of
$\sin(\pi z)/(\pi z)\in PW_\pi^2$, the space of square-integrable entire
functions of exponential type at most $\pi$. Likewise,
$\{\sqrt{4\pi k}:k\ge1\}$ is contained in the zero set of
$\sin(z^2/4)/(z^2/4)$ in the Fock space
$\{g\text{ entire}:\int_{\mathbb C}|g(z)|^2e^{-|z|^2}\,dA(z)<\infty\}$;
values at zero are defined by continuity. Lindel\"of's conditions require
bounded sums $\sum_{0<|a|\le R}a^{-1}$ for entire functions of exponential
type and $\sum_{0<|a|\le R}a^{-2}$ for functions of order two and finite
type, with zeros counted according to multiplicity. The positive integers
violate the first condition. In the Fock example, point evaluation gives
order at most two and finite type, while the sequence forces order at
least two and makes the second sum diverge. Positive multiplicities only
increase that sum.}
Our results establish the distinction in disk coefficient spaces,
including the analytic Wiener algebra, without fixing the multiplicities
in advance.

By Lemma~\ref{lem:division}, any finite collection of unwanted zeros of
a witness can be removed while remaining in $\ell_A^p$. Thus, if $W$ is
not exact, every nonzero function vanishing on $W$ must have infinitely
many zeros outside $W$. It is not enough to find one witness with extra
zeros: one must rule out every possible function whose zero set is exactly
$W$.

\subsection{Moment duality and exact polygon constructions}\label{subsec:intro-moment}
The moment-duality result below gives a sharp quantitative criterion for the vanishing problem. For $W\subset\D\setminus\{0\}$, let
\[
 A_p(W)=\inf_{f(0)=1,\ f|_W=0}\|f\|_p,
\]
where the infimum is over $f\in\ell_A^p$, with value $+\infty$ for an
empty class.

We establish a sharp duality formula for $A_p$ in terms of moments of
finitely supported complex measures on $W$. With the conjugate-exponent
convention above, for nonempty $W$ put
\[
 \delta_q(W)=\inf_{\widehat\mu(0)=1}
       \| (\widehat\mu(j))_{j\ge1}\|_{\ell^q},
 \qquad \widehat\mu(j)=\int_W w^j\,d\mu(w),
\]
where the infimum runs over such \emph{complex} measures. Theorem~\ref{thm:moment}
gives
\[
 W\text{ is a zero set for }\ell_A^p
 \quad\Longleftrightarrow\quad\delta_q(W)>0,
 \qquad A_p(W)^p=1+\delta_q(W)^{-p}.
\]
This includes $p=1$, where the proof uses Hahn--Banach on $c_0$, the
space of sequences tending to zero, with its supremum norm. The formula
is a dual version of the finite-extremal characterization in~\cite{CMR},
with the sharp normalized constant. It characterizes vanishing but does
not by itself guarantee exactness.

We next give exactness criteria for unions of regular polygons centered at the origin. For a natural family of these sets, we obtain a
necessary and sufficient summability condition for being a zero set, and
prove that the same condition also guarantees exact realization.
\begin{maintheorem}[An exact summability criterion]\label{main:rings}
Choose positive integers $d_m$, and define
\[
 n_1=1,\qquad n_{m+1}=(d_m+1)n_m,\qquad
 W(d)=\bigcup_{m\ge1}\{z\in\D:z^{n_m}=e^{-1/d_m}\}.
\]
For every $1<p<\infty$,
\[
 W(d)\text{ is a zero set for }\ell_A^p
 \quad\Longleftrightarrow\quad
 \sum_m d_m^{1-p}<\infty
 \quad\Longleftrightarrow\quad
 W(d)\text{ is an exact zero set for }\ell_A^p.
\]
One explicit infinite product realizes all admissible exponents
simultaneously, with simple zeros. Furthermore, $W(d)$ is Blaschke if
and only if $\sum_m d_m^{-1}<\infty$.
\end{maintheorem}

The proof and the coefficient-norm identity are given in
Theorem~\ref{thm:exactrings}. The products allow an arbitrary critical
exponent of summability, either included or excluded. They also provide
explicit non-Blaschke exact zero sets for $p>2$, Blaschke uniqueness
sets for $p<2$, and pairs of sequences with identical moduli but
different zero-set behavior. Random angular distributions yield further
Blaschke uniqueness sets. Compared with~\cite{Vin}, which establishes existence of non-Blaschke
zero sets for $p>2$, we give explicit exact realizations and a sharp
summability criterion.

\subsection{Unions and shift-invariant subspaces}
Finally, we investigate the behavior of zero sets under unions.
Non-exactness is a different phenomenon from failure of closure under
unions. For the latter we prove the following.
\begin{maintheorem}[Failure of closure under unions]\label{main:unions}
For every $p>2$, there is an exact zero set $W$ for $\ell_A^p$ such that
$W\cap(-W)=\varnothing$ and $W\cup(-W)$ is a uniqueness set for
$\ell_A^p$. Both exact witnesses may be chosen to have only simple zeros.
\end{maintheorem}

Theorem~\ref{main:unions} is proved in Theorem~\ref{union-thm:union-failure}.

At the endpoint $p=1$, the space $\ell_A^1$ is the classical analytic
Wiener algebra and a Banach algebra under multiplication. Hence the union
of two zero sets is again a zero set, in contrast with the failure of
closure under unions for $p>2$ proved here. The existence of non-exact
zero sets at $p=1$ is a separate issue: exact realization and closure
under finite unions are distinct properties.

The construction uses a Gaussian analytic function and its reflection.
Related Gaussian arguments for unions in Bergman and Fock spaces can be
found in~\cite{LZ}.

Selecting rotations of one such function gives a continuum of pairwise
disjoint exact zero sets whose pairwise unions are uniqueness sets
(Proposition~\ref{prop:continuum-zero-sets}). Their vanishing subspaces
yield a continuum of nonzero closed, infinite-dimensional right-shift-invariant
subspaces of $\ell^p(\mathbb N_0)$ with pairwise intersection $\{0\}$
(Corollary~\ref{cor:disjointshift}). This recovers Abakumov's existence result~\cite{Ab}.
Here $\mathbb N_0=\{0,1,2,\ldots\}$, and the right shift inserts a zero
as the first coordinate.

\subsection{Methods and organization}
The moment criterion follows from duality, and exact products on regular
polygons provide functions with controlled coefficient norms. Non-exactness
is obtained by selecting vertices according to quadratic characters,
which distinguish squares from non-squares modulo a prime. Quadratic
characters have also been used in the study of polynomial zeros, though in
different settings. For example, Fekete polynomials, whose coefficients are
given by Legendre symbols, have been studied with regard to the distribution
of their zeros~\cite{CGPS}.
In our construction, cancellation in the resulting finite Fourier
transforms makes the function small at omitted vertices; logarithmic
estimates then force additional zeros. At $p=1$, the configuration is
localized near one boundary point using a lacunary power series, whose
nonzero coefficients have rapidly increasing gaps. The union construction
instead uses Gaussian analytic functions.

Our applications of Clark-measure theory are confined to the elementary
case of the inner function $\vartheta(z)=z^n$, where the Clark transform
is a finite Fourier transform. We nevertheless include the general model-space viewpoint,
which readers familiar with the theory may find a convenient way to
organize several of the arguments. In this language, model-space
expansions explain folding of coefficients modulo $n$, while composition
explains the residue-class blocks of the character matrix. Averaging over the
Clark parameter, called disintegration, supplies the circle mean needed
for Jensen's formula. Section~\ref{sec:standard} develops these connections,
including the distinction between a fixed Clark parameter and almost-everywhere
averaging; the coefficient estimates are proved in Section~\ref{sec:lemmas}.

Sections~\ref{sec:obstructions}--\ref{sec:moduli} treat complete polygons
and angular dependence. The non-exactness constructions are in
Sections~\ref{sec:localnonexact}--\ref{sec:endpoint}, with the required
arithmetic introduced where it is used. Section~\ref{sec:unions} proves
the union result, and Section~\ref{sec:invariants} applies it to
shift-invariant subspaces.

{\bf AI statement.} The authors have been interested in questions considered in this paper for several years in connection with problems on the invariant subspaces of the shift operator in $\ell^p_A$, see for instance section \ref{sec:invariants}. The main ideas of statements and examples in the present text are human. LLM models, Anthropic's Fable 5 and ChatGPT's Astra, supplied the missing tricks and technical calculations to finalize our examples and theorems and helped create a latex presentation. The authors edited the latex file and verified the proofs to obtain the current version of the article.

\section{Analytic preliminaries and Clark measures}\label{sec:standard}
The general identities below provide context for the finite-dimensional
formulas used later. The discussion of the adjoint of composition in
Section~\ref{subsec:averaging-composition} and the Poisson identity at the
end of Section~\ref{subsec:logarithmic-averages} are included for perspective
and are not needed for the proofs of the main results.

\subsection{Sequence inequalities}
We write $D(c,R)=\{z:|z-c|<R\}$. A zero of order $m$ is a point
where the derivatives of orders $0,\ldots,m-1$ vanish and the derivative
of order $m$ does not.
We use counting-measure norms on coefficient vectors and normalized measure
on circles. The notation $a\asymp b$ means that their ratio is bounded
above and below by positive constants in the indicated limit. H\"older's inequality is
$\sum|a_jb_j|\le\|a\|_p\|b\|_q$, and Minkowski's inequality is
$\|\sum_i x^{(i)}\|_p\le\sum_i\|x^{(i)}\|_p$.
For $1\le p\le2$,
\[
 \sum_j|x_j|^2\le(\max_j|x_j|)^{2-p}\sum_j|x_j|^p,\qquad
 \|x\|_2\le\|x\|_p.
\]
For $x\in\mathbb C^n$, H\"older also gives
$\|x\|_p\le n^{1/p-1/2}\|x\|_2$. For probability weights $\lambda_j$,
\[
 \sum_j\lambda_j\log|Y_j|
 \le\frac12\log\sum_j\lambda_j|Y_j|^2,
\]
with $\log0=-\infty$. We use
\[
 1-r\le\log(1/r)\le(1-r)/r,\quad
 1-e^{-\varepsilon}\ge\varepsilon/2\ (0\le\varepsilon\le1).
\]

\subsection{Power series and removal of finitely many zeros}
H\"older gives, for $p>1$,
\begin{equation}\label{eq:evaluation}
 |f(z)|\le\|f\|_p(1-|z|^q)^{-1/q}.
\end{equation}
For $p=1$, $|f(z)|\le\|f\|_1$ on the closed disk, and $f$ is continuous there.
Every norm-bounded sequence in $\ell_A^p$ consequently has a subsequence
converging uniformly on compact subsets of $\D$, by Montel's theorem.
Coefficientwise convergence and Fatou's lemma show that the norm of the
limit is at most the lower limit of the norms.

Parseval's identity is
\[
 M_2(f,r)^2=\int_{\T}|f(r\zeta)|^2\dm(\zeta)
             =\sum_{j\ge0}|a_j|^2r^{2j}.
\]
Here and below $\mT$ denotes normalized Lebesgue measure on $\T$.
Nonzero $H^2$ functions have Blaschke zeros, and every Blaschke sequence has its usual bounded Blaschke product;
we use the normalization
\[
 b_a(z)=\frac{|a|}{a}\frac{a-z}{1-\bar az},\qquad
 B(0)=\prod_a|a|>0
\]
when zero is not prescribed. See, for example, \cite{Du,Nik} for these Hardy-space facts.

\begin{lemma}[Finite zero removal]\label{lem:division}
If $f\in\ell_A^p$, $1\le p<\infty$, and $f(a)=0$, $a\in\D$, then
\[
 \frac{f(z)}{z-a}\in\ell_A^p,\qquad
 \left\|\frac{f}{z-a}\right\|_p\le\frac{\|f\|_p}{1-|a|}.
\]
Consequently no point outside a zero set is a forced common zero of all its
vanishing witnesses. If a zero set is not exact, each nonzero witness has
infinitely many unwanted zeros.
\end{lemma}
\begin{proof}
Let $B_0f=(f-f(0))/z$ be the coefficient backward shift. Then
\[
 \frac{f(z)-f(a)}{z-a}
  =\sum_{k\ge0}a^k B_0^{k+1}f,
\]
with convergence in coefficient norm. Each shift has norm at most one,
which proves the estimate. Repeated division removes the full finite
multiplicity of any unwanted zero, without introducing new zeros.
\end{proof}

\subsection{Roots of unity and finite Fourier transforms}\label{subsec:dft}
Let $\omega=e^{2\pi i/n}$, and index vectors by $\mathbb Z_n=\{0,\ldots,n-1\}$.
The geometric sum gives
\[
 \sum_{k=0}^{n-1}\omega^{jk}
 =\begin{cases}n,&n\mid j,\\0,&n\nmid j.\end{cases}
\]
Consequently
\begin{equation}\label{eq:dft}
 Y_k=\sum_{l=0}^{n-1}c_l\omega^{lk},\qquad
 c_l=\frac1n\sum_{k=0}^{n-1}Y_k\omega^{-lk},\qquad
 \frac1n\sum_k|Y_k|^2=\sum_l|c_l|^2.
\end{equation}
These follow by substitution and expansion of the square. The unitary
transform is the first map divided by $\sqrt n$. Pointwise multiplication
of two sample vectors corresponds to cyclic convolution of their coefficient
vectors; hence its coefficient norm is bounded by the product of an
$\ell^1$ norm and an $\ell^p$ norm.

\subsection{Complex-analysis tools}
We use Rouch\'e's theorem, the sub-mean value property, Harnack's
inequality, and Jensen's formula. The needed normalizations are
\[
 |h(x)|^2\le\frac1{\pi d^2}\int_{D(x,d)}|h|^2\,dA,\qquad
 u(x)\ge\frac{R-|x-c|}{R+|x-c|}u(c)
\]
for nonnegative harmonic $u$ on $D(c,R)$. For analytic $h$, $\log|h|$
is subharmonic, and it is harmonic in a zero-free region.
For $h(0)\ne0$, Jensen's formula is
\[
 \sum_{0<|a|<r}\log\frac r{|a|}
  =\int_{\T}\log|h(r\zeta)|\dm(\zeta)-\log|h(0)|,
\]
with actual multiplicities. The logarithm on a fixed circle is integrable
even when the circle contains zeros.

We will sometimes use the analytic function
\begin{equation}\label{eq:orbitproduct}
 \mathcal N_h(x)=\prod_{z^n=x}h(z),\qquad \mathcal N_h(0)=h(0)^n.
\end{equation}
Away from zero take a local branch of $x^{1/n}$; changing the branch only
permutes the factors. Boundedness near zero makes its apparent singularity
removable. Thus $n^{-1}\log|\mathcal N_h|$ is subharmonic.

\subsection{Clark measures and model spaces}\label{subsec:clark}
Throughout this discussion, $\vartheta$ is an inner function with
$\vartheta(0)=0$. For $\alpha\in\T$, its Clark measure
$\sigma_\alpha^\vartheta$ is characterized by
\[
 \frac{\alpha+\vartheta(z)}{\alpha-\vartheta(z)}
 =\int_\T\frac{\zeta+z}{\zeta-z}\,d\sigma_\alpha^\vartheta(\zeta).
\]
We omit the superscript when the inner function is fixed. The normalization
at zero makes every $\sigma_\alpha$ a probability measure. It is singular
with respect to $\mT$ and is carried by the set where the nontangential
boundary value of $\vartheta$ is $\alpha$. Here a nontangential limit at
$\zeta\in\T$ means a limit as $z\to\zeta$ within
$|z-\zeta|\le C(1-|z|)$, for each fixed $C>1$.

The \emph{model space} is
$K_\vartheta=H^2\ominus\vartheta H^2$, the orthogonal complement of
$\vartheta H^2$ in $H^2$. Standard Clark-transform and boundary-value results identify it unitarily
with $L^2(\sigma_\alpha)$ by nontangential boundary restriction;
see~\cite{CMRb,Cl,GMR,Po,PS}.
We denote this restriction map by $V_\alpha$ and call it the
\emph{Clark transform}. In particular, for $g,h\in K_\vartheta$,
\begin{equation}\label{eq:clark-isometry-general}
 \int_\T (V_\alpha g)\overline{(V_\alpha h)}\,d\sigma_\alpha
       =\langle g,h\rangle_{H^2}.
\end{equation}
Our inner products are linear in the first argument. When integrating
members of $K_\vartheta$, we usually write $g$ instead of $V_\alpha g$.
The general identities below explain the constructions; their applications
in this paper use monomials $\vartheta(z)=z^n$, for which the transforms
are finite Fourier transforms.

\subsection{The measures and their moments}\label{subsec:clarkmoments}
The basic mean identity is
\begin{equation}\label{eq:clark-mean}
 \int_\T g\,d\sigma_\alpha=g(0),\qquad g\in K_\vartheta.
\end{equation}
Indeed $1\in K_\vartheta$, because
$\langle\vartheta h,1\rangle=(\vartheta h)(0)=0$ for every $h\in H^2$.
Also $V_\alpha1=1$, so~\eqref{eq:clark-isometry-general} with $h=1$
gives~\eqref{eq:clark-mean}. Integrability follows from Cauchy--Schwarz
and $\sigma_\alpha(\T)=1$. Thus a Clark measure represents evaluation
at zero on the model space, not on all of $H^2$.
More generally, if $\vartheta$ has a zero of order at least $n$ at zero,
then $1,z,\ldots,z^{n-1}\in K_\vartheta$, and consequently
\[
 \int_\T\zeta^j\,d\sigma_\alpha(\zeta)=0\qquad(1\le j<n).
\]

For $\vartheta(z)=z^n$, write $\sigma_\alpha^{(n)}$ for its Clark measure.
Then $K_{z^n}=\operatorname{span}\{1,z,\ldots,z^{n-1}\}$. Since
$\zeta^n=\alpha$ almost everywhere for $\sigma_\alpha^{(n)}$, writing
$j=kn+l$, $0\le l<n$, reduces its $j$th moment to $\alpha^k$ times its
$l$th moment. We obtain
\begin{equation}\label{eq:clarkmeasure}
 \sigma_\alpha^{(n)}=\frac1n\sum_{\zeta^n=\alpha}\delta_\zeta,\qquad
 \int_\T\zeta^j\,d\sigma_\alpha^{(n)}(\zeta)
 =\begin{cases}\alpha^{j/n},&n\mid j,\\0,&n\nmid j.\end{cases}
\end{equation}
Here $\delta_\zeta$ is unit point mass and the moment formula holds for
all integers $j$; negative moments follow by conjugation. The atomic
formula follows as well: the measure is carried by the $n$ roots of
$\alpha$, and its moments of orders $0,\ldots,n-1$ determine their masses
by finite Fourier inversion~\eqref{eq:dft}; these masses are all $1/n$.

For $G(z)=\sum_{l<n}c_lz^l$, finite Parseval~\eqref{eq:dft} gives
\[
 \int_\T|G|^2\,d\sigma_\alpha^{(n)}
       =\sum_{l<n}|c_l|^2=\|G\|_{H^2}^2.
\]
Thus finite Parseval is precisely the Clark isometry in this case.
Its inverse is
\[
 u\longmapsto(1-\bar\alpha z^n)
       \int_\T\frac{u(\zeta)}{1-\bar\zeta z}\,
                    d\sigma_\alpha^{(n)}(\zeta),
\]
whose coefficient of $z^l$, $0\le l<n$, is
$\int u(\zeta)\bar\zeta^l\,d\sigma_\alpha^{(n)}(\zeta)$.
A polygon inside the disk is not itself a Clark measure on $\T$.
Its normalized counting measure is the \emph{radial image} of
$\sigma_\alpha^{(n)}$: move each atom $\zeta$ to $r\zeta$ without
changing its mass. This observation will give the moment tests in
Section~\ref{sec:obstructions}.

\subsection{Model-space expansions and folding}\label{subsec:clarkfold}
Every $f\in H^2$ has a unique orthogonal expansion
\begin{equation}\label{eq:theta-expansion}
 H^2=\bigoplus_{k\ge0}\vartheta^kK_\vartheta,\qquad
 f=\sum_{k\ge0}\vartheta^kf_{k+1},\qquad
 \|f\|_{H^2}^2=\sum_{k\ge0}\|f_{k+1}\|_{H^2}^2,
\end{equation}
where $f_{k+1}\in K_\vartheta$ and the series converges in $H^2$.
To see this, iterate $H^2=K_\vartheta\oplus\vartheta H^2$.
The remaining intersection $\bigcap_{k\ge0}\vartheta^kH^2$ is zero,
since any member would vanish to arbitrarily high order at zero.
This proves completeness of the orthogonal sum.

For a finite sum $f=\sum_{k=0}^M\vartheta^kf_{k+1}$, the identity
$\vartheta=\alpha$ $\sigma_\alpha$-almost everywhere gives
\begin{equation}\label{eq:finite-clarkfold}
 f=V_\alpha\left(\sum_{k=0}^M\alpha^kf_{k+1}\right)
       \quad\sigma_\alpha\text{-a.e.},\qquad
 \int_\T f\,d\sigma_\alpha=\sum_{k=0}^M\alpha^kf_{k+1}(0).
\end{equation}
The function $\sum_k\alpha^kf_{k+1}\in K_\vartheta$ is its
\emph{fold at $\alpha$}: on the Clark level set the powers of
$\vartheta$ are replaced by the scalars $\alpha^k$.
For an infinite expansion the same identities hold whenever its partial
sums converge to the boundary function $f$ in $L^2(\sigma_\alpha)$.
One sufficient condition, valid for every $\alpha$, is
$\sum_{k\ge0}\|f_{k+1}\|_{H^\infty}<\infty$.
Indeed the expansion is then uniformly convergent in the disk, and its
boundary series is absolutely convergent off the countable union of the
exceptional sets for its model-space terms. Under this condition,
\begin{equation}\label{eq:clarkfold-mean}
 \int_\T f\,d\sigma_\alpha=\sum_{k\ge0}\alpha^kf_{k+1}(0),\qquad \alpha\in\T.
\end{equation}

The distinction between $H^2$ convergence and convergence for a fixed
Clark measure is necessary. For $\vartheta(z)=z$ one has
$K_\vartheta=\mathbb C$ and $\sigma_1=\delta_1$. The function
$f(z)=\sum_{k\ge0}z^k/(k+1)$ belongs to $H^2$, but $f(r)\to+\infty$
as $r\uparrow1$ and the corresponding sum in~\eqref{eq:clarkfold-mean}
diverges. Thus the formula at the prescribed parameter $\alpha=1$ is
not an unconditional identity for all $H^2$ functions.
Section~\ref{subsec:averaging-composition} gives its general almost-everywhere
interpretation.

In every application below this boundary issue is avoided by dilating
first. Let $f(z)=\sum_{j\ge0}a_jz^j$ be analytic in $\D$ and $0<r<1$.
Then $f_r(z)=f(rz)$ has absolutely summable coefficients and, for
$\vartheta(z)=z^n$, its blocks are
\[
 g_k(z)=\sum_{l=0}^{n-1}a_{kn+l}r^{kn+l}z^l\in K_{z^n},\qquad
 f_r=\sum_{k\ge0}z^{kn}g_k,\qquad
 \sum_k\|g_k\|_{H^\infty}\le\sum_j|a_j|r^j<\infty.
\]
Its fold is therefore the polynomial
\begin{equation}\label{eq:monomial-fold}
 G_{r,\alpha}(z)=\sum_{k\ge0}\alpha^kg_k(z)
  =\sum_{l=0}^{n-1}c_lz^l,\qquad
 c_l=\sum_{k\ge0}\alpha^ka_{kn+l}r^{kn+l}.
\end{equation}
It agrees with $f_r$ at every atom of $\sigma_\alpha^{(n)}$, and
\begin{equation}\label{eq:fold-pairing}
 \int_\T f_r\overline h\,d\sigma_\alpha^{(n)}
     =\langle G_{r,\alpha},h\rangle_{H^2}\quad(h\in K_{z^n}),\qquad
 \int_\T|f_r|^2\,d\sigma_\alpha^{(n)}=\sum_{l<n}|c_l|^2.
\end{equation}
In particular the mean is $c_0=\sum_{k\ge0}\alpha^ka_{kn}r^{kn}$;
testing against $h(z)=z^l$ recovers every folded coefficient $c_l$.
The scalar mean identity is thus the constant coefficient of the
full folding identity, and the Clark isometry supplies its mean-square
version. Lemma~\ref{lem:fold} adds the coefficient-norm bounds.

\subsection{Spectral projections and composition}\label{subsec:clarkspectral}
For the monomial model space the Clark unitary is
\[
 U_\alpha=V_\alpha^{-1}M_\zeta V_\alpha,\qquad
 U_\alpha z^l=z^{l+1}\ (l<n-1),\qquad U_\alpha z^{n-1}=\alpha,
\]
where $M_\zeta$ is multiplication by $\zeta$ on
$L^2(\sigma_\alpha^{(n)})$. For a function $\phi$ on the atoms, define
$\phi(U_\alpha)=V_\alpha^{-1}M_\phi V_\alpha$. Then
\begin{equation}\label{eq:clarkmatrix}
 \langle\phi(U_\alpha)z^{l'},z^l\rangle
   =\int_\T\phi(\zeta)\zeta^{l'-l}\,d\sigma_\alpha^{(n)}(\zeta).
\end{equation}
In particular, selecting a subset of atoms gives an orthogonal
\emph{spectral projection}, the conjugate of multiplication by its
indicator. The moments of that indicator are the matrix entries of this
projection. Gauss-sum cancellation will make its off-diagonal entries
small, preventing a vector in its range from concentrating on too few
monomials. This is a Hilbert-space statement; it does not identify
coefficient $\ell^p$ norms with Clark $L^p$ norms when $p\ne2$.

The grouping of atoms for composite polygon sizes also has a general
model-space interpretation. For every integer $N\ge1$,
\begin{equation}\label{eq:clark-power-composition}
 K_{\vartheta^N}=\bigoplus_{k=0}^{N-1}\vartheta^kK_\vartheta,
 \qquad
 \sigma_1^{\vartheta^N}=\frac1N\sum_{\beta^N=1}\sigma_\beta^\vartheta.
\end{equation}
The first identity follows by iterating the model-space decomposition.
For the second use the geometric-series identity
\[
 \frac{1+w^N}{1-w^N}
    =\frac1N\sum_{\beta^N=1}\frac{\beta+w}{\beta-w},\qquad |w|<1,
\]
put $w=\vartheta(z)$, and compare Herglotz representations.
Taking $\vartheta(z)=z^{n'}$ and $n=n'N$ yields
\begin{equation}\label{eq:clarkcomposition}
 \sigma_1^{(n)}=\frac1N\sum_{\beta^N=1}\sigma_\beta^{(n')}.
\end{equation}
Thus the large polygon splits into the level sets of $z^{n'}$.
Combining this identity with the zero moments on $K_{z^{n'}}$ will
produce the residue-class blocks in the character matrix of
Lemma~\ref{lem:character}.

\subsection{Disintegration}\label{subsec:disintegration}
Averaging the Clark measures over their parameter recovers normalized
Lebesgue measure on the circle. This identity is called
\emph{disintegration}. For any inner $\vartheta$ with $\vartheta(0)=0$
and $u\in L^1(\mT)$, with a measurable representative, it says
\begin{equation}\label{eq:disintegration}
 \int_\T\left(\int_\T u\,d\sigma_\alpha^\vartheta\right)\dm(\alpha)
      =\int_\T u\dm.
\end{equation}
The inner integral is defined for almost every $\alpha$; applying the
identity to $|u|$ gives the required integrability. This is Aleksandrov's
disintegration theorem~\cite{Al,PS}. Under our normalization its
measure identity also follows by averaging the Herglotz representations:
for each fixed $z$, the average of
$(\alpha+\vartheta(z))/(\alpha-\vartheta(z))$ is $1$. Uniqueness of the
representing measure gives $\int\sigma_\alpha^\vartheta\dm(\alpha)=\mT$;
positivity then extends the formula from continuous functions to $L^1$.
For monomials it is simply a change of variables:
\begin{align*}
 \int_\T\left(\int_\T u\,d\sigma_\alpha^{(n)}\right)\dm(\alpha)
 &=\frac1{2\pi n}\sum_{k=0}^{n-1}\int_0^{2\pi}
       u\bigl(e^{i(t+2\pi k)/n}\bigr)\,dt\\
 &=\frac1{2\pi}\sum_{k=0}^{n-1}
       \int_{2\pi k/n}^{2\pi(k+1)/n}u(e^{is})\,ds
 =\int_\T u\dm.
\end{align*}

\subsection{Averaging as the adjoint of composition}\label{subsec:averaging-composition}
For $f\in H^2$, use its boundary values and set
$\mathcal A_\vartheta f(\alpha)=\int_\T f\,d\sigma_\alpha^\vartheta$
for almost every $\alpha$. Cauchy--Schwarz and disintegration give
$\|\mathcal A_\vartheta f\|_{L^2(\mT)}\le\|f\|_{H^2}$.
For $j\ge0$, the support property $\vartheta=\alpha$
$\sigma_\alpha^\vartheta$-almost everywhere gives
\[
 \int_\T\mathcal A_\vartheta f(\alpha)\bar\alpha^j\dm(\alpha)
 =\int_\T f\,\overline{\vartheta^j}\dm
 =\langle f,\vartheta^j\rangle_{H^2}=f_{j+1}(0),
\]
where the last equality uses~\eqref{eq:theta-expansion} and
$1\in K_\vartheta$. Its negative Fourier coefficients vanish, since
$\int f\vartheta^j\dm=(f\vartheta^j)(0)=0$ for $j\ge1$.
Consequently, with equality in $L^2(\mT)$,
\begin{equation}\label{eq:clark-average-H2}
 \mathcal A_\vartheta f(\alpha)
    =\sum_{j\ge0}\alpha^jf_{j+1}(0),\qquad
 \mathcal A_\vartheta|_{H^2}=C_\vartheta^*,\qquad
 C_\vartheta h=h\circ\vartheta.
\end{equation}
The composition operator $C_\vartheta$ is an isometry because
$1,\vartheta,\vartheta^2,\ldots$ are orthonormal. The coefficient
calculation proves the adjoint formula directly; this is the
Aleksandrov averaging operator, as in~\cite{GP,PS}.
Pointwise, the series in~\eqref{eq:clark-average-H2} can always be
interpreted by its Abel limit
$\lim_{s\uparrow1}\sum_{j\ge0}s^j\alpha^jf_{j+1}(0)$ for almost every
$\alpha$. This is the general $H^2$ version of the mean formula; it does
not assert convergence at a prescribed parameter such as $\alpha=1$.
For $\vartheta(z)=z^n$ it reads
$C_{z^n}^*f(w)=\sum_{j\ge0}a_{nj}w^j$: averaging keeps exactly the
constant coefficients of the model-space blocks.

\subsection{Logarithmic averages}\label{subsec:logarithmic-averages}
For $0<\rho<1$ and analytic $h\not\equiv0$, let
\[
 L_h(\alpha)=\frac1n\sum_{z^n=\rho\alpha}\log|h(z)|
  =\int_\T\log|h(\rho^{1/n}\zeta)|\,d\sigma_\alpha^{(n)}(\zeta).
\]
The integrand is integrable with respect to $\mT$, even if its circle
contains a zero. Disintegration and Jensen's formula give
\begin{equation}\label{eq:clark-log-average}
 \int_\T L_h\dm
 =\int_\T\log|h(\rho^{1/n}\zeta)|\dm(\zeta)
 \ge\log|h(0)|.
\end{equation}
Equivalently, $L_h(\alpha)=n^{-1}\log|\mathcal N_h(\rho\alpha)|$.
Unlike~\eqref{eq:clark-mean}, this argument applies to the nonanalytic
function $\log|h|$: averaging over the Clark parameter, rather than
integration against a single Clark measure, is what permits comparison
with Jensen's formula. This distinction is used in the proofs of
non-exactness in Sections~\ref{sec:localnonexact} and~\ref{sec:endpoint}.

As a further illustration, the Poisson kernel
$P_r(t)=(1-r^2)/|1-re^{it}|^2$, $0\le r<1$, satisfies
\[
 \frac1n\sum_{j=0}^{n-1}P_r(t-2\pi j/n)=P_{r^n}(nt).
\]
Indeed expand $P_r(t)=1+2\sum_{j\ge1}r^j\cos jt$ and
apply~\eqref{eq:clarkmeasure}. Thus the sparse moments, coefficient
folding, and Poisson averaging are three consequences of the same
Clark-measure identities.

\section{Moment duality and coefficient constructions}\label{sec:lemmas}
This section establishes a dual criterion for zero sets and develops the
coefficient estimates used in the explicit constructions that follow.

\subsection{Moment duality, including the best constant}
Recall from Section~\ref{subsec:intro-moment} that
$\widehat\mu(j)=\int w^j\,d\mu(w)$ for a finitely supported complex
measure $\mu$ on $W\subset\D\setminus\{0\}$, and that, for nonempty $W$,
\[
 \delta_q(W)=\inf_{\widehat\mu(0)=1}
       \|(\widehat\mu(j))_{j\ge1}\|_q.
\]
For $q=\infty$ the norm is the supremum over $j\ge1$. We use the
conventions $\delta_q(\varnothing)=+\infty$ and
$A_p(\varnothing)=1$.
The infimum is over complex measures. Positive measures suffice for several
obstructions below but do not replace complex measures in this equivalence.

\begin{theorem}[The sharp moment theorem]\label{thm:moment}
Let $1\le p<\infty$. A nonempty set $W\subset\D\setminus\{0\}$ is a zero set
for $\ell_A^p$ if and only if $\delta_q(W)>0$. In that case
\begin{align}
 \min_{f(0)=1,\ f|_W=0}\|f-1\|_p&=\delta_q(W)^{-1},\label{eq:tailmoment}\\
 A_p(W)=\min_{f(0)=1,\ f|_W=0}\|f\|_p
    &=\bigl(1+\delta_q(W)^{-p}\bigr)^{1/p}.\label{eq:fullmoment}
\end{align}
For $p>1$ the minimizer is unique. No uniqueness assertion is made at $p=1$.
\end{theorem}
\begin{proof}
Any nonzero vanishing function can be divided by its first nonzero monomial
and normalized to have constant coefficient one, since $0\notin W$.
If $f=1+\sum_{j\ge1}a_jz^j$ vanishes on $W$, then
\[
 1=\left|\sum_{j\ge1}a_j\widehat\mu(j)\right|
   \le\|f-1\|_p\|(\widehat\mu(j))_{j\ge1}\|_q.
\]
This proves necessity and the lower bound.

For the converse put $v_w=(w,w^2,\ldots)$ and, on their algebraic span, define
\[
 L\left(\sum_wc_wv_w\right)=-\sum_wc_w.
\]
The definition of $\delta_q$ says exactly that
$|\sum c_w|\le\delta_q(W)^{-1}\|\sum c_wv_w\|_q$, so $L$ is well defined,
bounded, and has norm exactly $\delta_q(W)^{-1}$.
For $p>1$ extend it by Hahn--Banach to $\ell^q$, whose dual is $\ell^p$.
For $p=1$, every $v_w$ lies in $c_0$, so instead extend to $c_0$ with its
supremum norm and use $(c_0)^*=\ell^1$. Extending to $\ell^\infty$ and
identifying its whole dual with $\ell^1$ would be incorrect.
The resulting coefficient sequence $(a_j)_{j\ge1}$ has the required norm
and satisfies $\sum_{j\ge1}a_jw^j=-1$. This proves
\eqref{eq:tailmoment}; adding the constant coefficient proves
\eqref{eq:fullmoment}. Strict convexity gives uniqueness for $p>1$.
\end{proof}

In particular,
\begin{equation}\label{eq:finitecriterion}
 W\text{ is a zero set}\quad\Longleftrightarrow\quad
       \sup_{\substack{E\subset W\\ E\ {\rm finite}}}A_p(E)<\infty.
\end{equation}
This follows because $\delta_q(W)=\inf_{E\subset W,\ E\ {\rm finite}}
\delta_q(E)$. The point zero can be inserted by multiplication by $z$.
For finite $E$ with distinct points and $p=2$, the normalized Blaschke product
is extremal and
\[
 A_2(E)=\prod_{w\in E}|w|^{-1},\qquad
 \delta_2(E)=\left(\prod_{w\in E}|w|^{-2}-1\right)^{-1/2}.
\]
For a regular polygon the positive test measure is a radial image of
$\sigma_1^{(n)}$. For a general $W$, however, the complex measures in the
duality theorem need not be Clark measures or positive measures.

\subsection{A convolution estimate}
\begin{lemma}\label{lem:convolution}
For $q\ge1$, $0<\gamma<1$, and nonnegative $x_1,\ldots,x_M$,
\[
 \sum_{v=1}^M\left(\sum_{m\le v}\gamma^{v-m}x_m\right)^q
 \le(1-\gamma)^{-q}\sum_{m=1}^Mx_m^q.
\]
\end{lemma}
\begin{proof}
For $q>1$, H\"older bounds the inner sum to the power $q$ by
$(\sum_{m\le v}\gamma^{v-m})^{q-1}
 \sum_{m\le v}\gamma^{v-m}x_m^q$.
Sum over $v$ and use the geometric series twice. For $q=1$, interchange
the two finite sums directly.
\end{proof}

\subsection{Products with disjoint coefficient contributions}
\begin{lemma}\label{lem:products}
Let $P_m(x)=\sum_{k=0}^{D_m}p_{m,k}x^k$, $P_m(0)=1$, $D_m\ge1$.
Suppose the positive integers $n_m$ satisfy
\begin{equation}\label{eq:superincreasing}
 n_m>\sum_{j<m}D_jn_j.
\end{equation}
For $F_M=\prod_{m\le M}P_m(z^{n_m})$ one has
\[
 \|F_M\|_p^p=\prod_{m\le M}\|P_m\|_p^p,\qquad
 \deg F_M<n_{M+1}.
\]
The coefficients stabilize; if $\sup_M\|F_M\|_p\le C$, the limiting power
series $F$ belongs to $\ell_A^p$, has $F(0)=1$, and the $F_M$ are Taylor
truncations of $F$. They converge locally uniformly, and all zeros of their
factors persist.
\end{lemma}
\begin{proof}
Every coefficient contribution has exponent $\sum k_mn_m$, $0\le k_m\le D_m$.
If two digit choices differ, their largest differing index contributes at
least $n_m$, whereas the sum of all preceding differences is smaller than
$n_m$ by \eqref{eq:superincreasing}. Thus no contributions overlap.
Summing the $p$-th powers factors the norm exactly.
All later factors start beyond $\deg F_M$, so coefficients stabilize and
$F_M$ is the appropriate truncation. Fatou gives $\sum|F_j|^p\le C^p$.
The estimate $\sum_{j\ge n_{M+1}}|F_j|r^j\le Cr^{n_{M+1}}/(1-r)$
gives local uniform convergence and persistence of the zeros.
\end{proof}
The sufficient condition $n_{m+1}\ge(D_m+1)n_m$
implies \eqref{eq:superincreasing}, since
$D_mn_m\le n_{m+1}-n_m$ and the sum telescopes. 

The polynomial $F_{M-1}$ lies in $K_{z^{n_M}}$.
In the expansion~\eqref{eq:theta-expansion} with
$\vartheta(z)=z^{n_M}$, the model-space blocks of
$F_M=P_M(z^{n_M})F_{M-1}$ are exactly $p_{M,k}F_{M-1}$.
At $p=2$, orthogonality therefore gives the norm factorization.
For general $p$, the same identity follows from the disjoint Taylor
coefficient blocks, not from a Clark-space isometry for $\ell^p$.

\subsection{Folding}
\begin{lemma}\label{lem:fold}
Let $n\ge1$ be an integer, let $f(z)=\sum_{j\ge0}a_jz^j$ be analytic
in $\D$, and let $0<|\xi|<1$. Set $\omega=e^{2\pi i/n}$ and
$t=|\xi|^n$. Put
\[
 C_l=\sum_{i\ge0}a_{l+in}\xi^{l+in}\quad(0\le l<n).
\]
Then
\[
 f(\xi\omega^k)=\sum_{l<n}C_l\omega^{lk},\qquad
 \frac1n\sum_{k<n}|f(\xi\omega^k)|^2=\sum_{l<n}|C_l|^2.
\]
If, in addition, $f\in\ell_A^p$ for $1\le p<\infty$, then
\begin{equation}\label{eq:fold-norm}
 \|C\|_p\le\kappa_p(t)\|f\|_p,\qquad
 \kappa_p(t)=
 \begin{cases}
  (1-t^q)^{-1/q},&1<p<\infty,\\
  1,&p=1.
 \end{cases}
\end{equation}
\end{lemma}
\begin{proof}
Write $\xi=r\beta$, $|\beta|=1$, and put $\alpha=\beta^n$.
In~\eqref{eq:monomial-fold}, the coefficients of the fold of $f(rz)$
satisfy $C_l=\beta^lc_l$. Its agreement with $f(rz)$ on the Clark atoms
and the isometry~\eqref{eq:fold-pairing} give the two sample identities.
The required absolute convergence was established in
Section~\ref{subsec:clarkfold}.

For the coefficient bound, $|\xi|^{l+in}\le t^i$ and H\"older's
inequality give, when $p>1$,
\[
 |C_l|\le\sum_{i\ge0}|a_{l+in}|t^i
 \le(1-t^q)^{-1/q}
       \left(\sum_{i\ge0}|a_{l+in}|^p\right)^{1/p}.
\]
Raise to the power $p$ and sum in $l$; the residue classes partition the
coefficient indices. For $p=1$, simply use $t^i\le1$ and sum the first
inequality. This proves~\eqref{eq:fold-norm}; the displayed $\kappa_p$
are finite and nondecreasing for $0\le t<1$.
\end{proof}

\section{Zero-set criteria and exact products for regular polygons}\label{sec:obstructions}
We combine a quantitative obstruction based on polygonal mass with explicit
products that realize sufficiently separated unions of regular polygons
exactly.

\subsection{Polygon mass and a quantitative obstruction}
Let
\[
 R_m=\{z:z^{n_m}=e^{-\varepsilon_m}\},\qquad
 r_m=e^{-\varepsilon_m/n_m},\qquad \varepsilon_m>0.
\]
Thus $R_m$ has $n_m$ equally spaced vertices. If $\varepsilon_m\le1$, its
Blaschke contribution lies between $\varepsilon_m/e$ and $\varepsilon_m$:
use $r_m\ge e^{-1}$ and
$r_m\log(1/r_m)\le1-r_m\le\log(1/r_m)$.
For its normalized counting measure $\nu_m$,
\begin{equation}\label{eq:ringmoments}
 \widehat\nu_m(j)=\boldsymbol1_{\{n_m\mid j\}}
                    e^{-\varepsilon_mj/n_m}.
\end{equation}

This is~\eqref{eq:clarkmeasure} for the radial image of
$\sigma_1^{(n_m)}$. Equivalently, the block-mean identity from
Section~\ref{subsec:clarkfold} gives, for analytic $f(z)=\sum_ja_jz^j$,
\[
 \int f\,d\nu_m=\sum_{k\ge0}a_{kn_m}r_m^{kn_m},
\]
the sum of the constant terms of the dilated model-space blocks.
A normalized function vanishing on $R_m$ makes this sum zero.
The moment obstruction below combines these block-mean constraints
across the different levels.

\begin{proposition}[Polygons that are too heavy]\label{prop:heavy}
Assume $n_m\mid n_{m+1}$ and $n_{m+1}\ge2n_m$.
For $1<p<\infty$, put $S_M=\sum_{m\le M}\varepsilon_m^{p-1}$.
Then
\[
 \delta_q\left(\bigcup_{m\le M}R_m\right)
       \le\frac{q^{-1/q}}{1-2^{-1/q}}S_M^{-1/p}.
\]
Consequently, every normalized vanishing witness satisfies
\[
 \|f-1\|_p\ge q^{1/q}(1-2^{-1/q})S_M^{1/p}.
\]
In particular, $\sum_m\varepsilon_m^{p-1}=\infty$ makes the union a uniqueness set.
\end{proposition}
\begin{proof}
Choose nonnegative $\alpha_m$ of sum one and put $\mu=\sum_{m\le M}\alpha_m\nu_m$.
Each nonzero moment can be grouped by the largest $v\le M$ for which
$n_v\mid j$, and then $j=n_vt$. Enlarging the allowed positive integers $t$,
\begin{align*}
 \sum_{j\ge1}|\widehat\mu(j)|^q
 &\le\sum_{v=1}^M\sum_{t\ge1}
       \left(\sum_{m\le v}\alpha_m
              e^{-\varepsilon_m(n_v/n_m)t}\right)^q\\
 &\le\frac1q\sum_{v=1}^M
       \left(\sum_{m\le v}2^{-(v-m)/q}
                     \alpha_m\varepsilon_m^{-1/q}\right)^q\\
 &\le\frac1{q(1-2^{-1/q})^q}
                 \sum_{m\le M}\frac{\alpha_m^q}{\varepsilon_m}.
\end{align*}
The second line is Minkowski in $t$ followed by
$\sum_{t\ge1}e^{-at}\le1/a$ and $n_v/n_m\ge2^{v-m}$.
The third is Lemma \ref{lem:convolution}.
Take $\alpha_m=\varepsilon_m^{p-1}/S_M$; since $(p-1)q=p$, the last sum is
$S_M^{1-q}$. Take $q$-th roots and apply Theorem \ref{thm:moment}.
\end{proof}

\subsection{Exact products for separated polygons}
The following factors are chosen so that their remaining polynomial
factors have no zeros in the disk. Thus the infinite product prescribes
the complete zero set, not merely a subset of it.

\begin{proposition}[Exact witnesses for sufficiently separated light polygons]
\label{prop:light}
Let $0<\varepsilon_m\le1$, $D_m=\lceil1/\varepsilon_m\rceil$,
$\rho_m=e^{-\varepsilon_m}$, and assume \eqref{eq:superincreasing}.
If $\sum_m\varepsilon_m^{p-1}<\infty$, $1<p<\infty$, then the single product
\[
 F(z)=\prod_{m\ge1}Q_m(z^{n_m}),\qquad
 Q_m(x)=1-\frac{x+x^2+\cdots+x^{D_m}}{S_m},\quad
 S_m=\sum_{k=1}^{D_m}\rho_m^k,
\]
has exactly the simple zeros $\bigcup_mR_m$ and belongs to $\ell_A^p$.
More precisely,
\[
 \|F\|_p^p=\prod_m\left(1+\frac{D_m}{S_m^p}\right)
       \le\exp\left(2\cdot5^p\sum_m\varepsilon_m^{p-1}\right).
\]
The function does not depend on the admissible exponent $p$.
\end{proposition}
\begin{proof}
Since $\varepsilon_mD_m\ge1$, $\rho_m\ge e^{-1}$, and $1-\rho_m\le\varepsilon_m$,
\[
 S_m=\rho_m\frac{1-\rho_m^{D_m}}{1-\rho_m}
       \ge\frac{e^{-1}(1-e^{-1})}{\varepsilon_m}
       >\frac1{5\varepsilon_m}.
\]
Also $D_m\le2/\varepsilon_m$. Thus $D_m/S_m^p\le2\cdot5^p\varepsilon_m^{p-1}$,
and Lemma \ref{lem:products} gives the norm identity and convergence.
For $|z|\le R<1$,
\[
 |Q_m(z^{n_m})-1|\le5\,\frac{R^{n_m}}{1-R^{n_m}}.
\]
Condition \eqref{eq:superincreasing} implies exponential growth of $n_m$,
so these deviations are locally summable.

Suppressing $m$, factor
\begin{equation}\label{eq:Qfactor}
 Q(x)=(1-x/\rho)B(x),\qquad
 B(x)=\sum_{j=0}^{D-1}b_jx^j,\quad
 b_j=\frac{\rho+\cdots+\rho^{D-j}}S.
\end{equation}
Here $b_0=1>b_1>\cdots>b_{D-1}>0$ for $D>1$, and
\[
 (1-x)B(x)=b_0-\sum_{j=1}^{D-1}(b_{j-1}-b_j)x^j-b_{D-1}x^D.
\]
The subtracted terms have total absolute value strictly less than $b_0$
for $|x|<1$. Thus $B$ is zero-free there, and $Q$ has only its simple zero
$\rho$. The case $D=1$ has $B=1$.
The polygon radii are distinct: $n_{m+1}>D_mn_m\ge n_m/\varepsilon_m$ gives
$\varepsilon_{m+1}/n_{m+1}<\varepsilon_m/n_m$.
The locally summable deviations exclude any extra zeros of the infinite
product. This proves exactness and simplicity.
\end{proof}

\begin{theorem}[An exact summability criterion]\label{thm:exactrings}
Choose positive integers $d_m$ and put
\[
 n_1=1,\qquad n_{m+1}=(d_m+1)n_m,\qquad
 W(d)=\bigcup_{m\ge1}\{z:z^{n_m}=e^{-1/d_m}\}.
\]
For $1<p<\infty$ the following are equivalent:
\[
 W(d)\text{ is a zero set for }\ell_A^p
 \quad\Longleftrightarrow\quad
 \sum_md_m^{1-p}<\infty
 \quad\Longleftrightarrow\quad
 W(d)\text{ is an exact zero set for }\ell_A^p.
\]
One product realizes all the admissible exponents, with simple zeros, and
\[
 \|F\|_p^p=
 \prod_m\left(1+\frac{d_m}{(\sum_{k=1}^{d_m}e^{-k/d_m})^p}\right)
 \le\exp\left(e^p\sum_md_m^{1-p}\right).
\]
Moreover $W(d)$ is Blaschke if and only if $\sum_md_m^{-1}<\infty$.
\end{theorem}
\begin{proof}
Necessity is Proposition \ref{prop:heavy} with $\varepsilon_m=1/d_m$.
For sufficiency use Proposition \ref{prop:light}, noting now
$d_m/e\le S_m\le d_m$. The recurrence makes \eqref{eq:superincreasing}
hold by telescoping. Finally
\[
 \frac{1-e^{-1}}{d_m}
 \le n_m(1-e^{-1/(d_mn_m)})\le\frac1{d_m},
\]
which gives the Blaschke assertion.
\end{proof}
The stronger conclusion is special to this explicitly separated family.
Proposition \ref{prop:heavy} needs only divisibility and doubling; the
arguments here do not settle sufficiency without the degree-separation
hypothesis, for instance for all dyadic polygon sizes.

\subsection{A Blaschke uniqueness set for every \texorpdfstring{$p<2$}{p < 2}}
The Blaschke condition does not characterize zero sets for $1<p<2$~\cite{CD}.
The following explicit polygonal construction works simultaneously for all
$p<2$.
\begin{corollary}\label{thm:dyadic}
Let
\[
 n_m=2^m,\qquad \varepsilon_m=\frac1{m\log^2(m+2)},\qquad
 W=\bigcup_m\{z:z^{2^m}=e^{-\varepsilon_m}\}.
\]
Then $W$ is Blaschke, is a uniqueness set for every $\ell_A^p$, $1\le p<2$,
and is exact for every $p\ge2$.
\end{corollary}
\begin{proof}
The Blaschke sum is bounded by $\sum_m\varepsilon_m<\infty$.
For $1<p<2$,
$m^{2-p}\ge\log^2(m+2)$ eventually, and
$\log^{2(p-1)}(m+2)\le\log^2(m+2)$, so
$\varepsilon_m^{p-1}\ge1/m$ eventually.
Proposition \ref{prop:heavy} applies. The case $p=1$ follows by
$\ell_A^1\subset\ell_A^{3/2}$, and $p\ge2$ by the Blaschke product.
\end{proof}

\subsection{Arbitrary critical exponents}
\begin{corollary}[Prescribed critical exponent]\label{thm:critical}
Let $p_0>1$ and $a=1/(p_0-1)$. For $d_m=\lceil m^a\rceil$, the set
$W(d)$ of Theorem~\ref{thm:exactrings} is a zero set for $\ell_A^p$,
$1<p<\infty$, exactly when $p>p_0$, and it is exact for every such $p$.
For
\[
 d_m=\left\lceil m^a\log^b(m+2)\right\rceil,
 \qquad b>\frac1{p_0-1},
\]
the corresponding range is exactly $p\ge p_0$. In each case one product
realizes all the admissible exponents, with simple zeros.
\end{corollary}
\begin{proof}
By Theorem~\ref{thm:exactrings}, the criterion is
$\sum_m d_m^{1-p}<\infty$. In the first case this is the numerical series
$\sum_m m^{-a(p-1)}$; in the second it is comparable to
\[
 \sum_m m^{-a(p-1)}\log^{-b(p-1)}(m+2).
\]
At $p=p_0$ the power exponent is one and the logarithmic exponent is
$b(p_0-1)>1$. These comparisons give both assertions.
\end{proof}
The corollary gives explicit exact realizations of critical exponents;
nesting of the zero-set classes as the exponent varies was established in~\cite{C}.

\subsection{Random arguments}
Let $\varphi_k$ be independent random variables, each uniformly distributed
on $[0,2\pi)$. We write $\E$ for expectation and $\mathbb P$ for probability;
``almost surely'' means outside an event of probability zero.
For an integer $h\ge1$, expansion and independence give
\[
 \E\left|\sum_kc_ke^{i\varphi_k}\right|^{2h}
      \le h!\left(\sum_k|c_k|^2\right)^h.
\]
Only matching multiplicities survive; for each ordered $h$-tuple there are
at most $h!$ matching reorderings. Thus, for $2\le q\le2h$,
\begin{equation}\label{eq:steinhaus}
 \E\left|\sum_kc_ke^{i\varphi_k}\right|^q
 \le(h!)^{q/(2h)}\left(\sum_k|c_k|^2\right)^{q/2}.
\end{equation}

\begin{proposition}\label{thm:random}
Let $1<p\le2$ and $0<s_k<1$. If
$\sum_k(1-s_k)^{2/q}=\infty$, then almost surely
$\{s_ke^{i\varphi_k}\}$ is a uniqueness set for $\ell_A^p$.
\end{proposition}
\begin{proof}
Put $t_k=1-s_k$, $T_N=\sum_{k\le N}t_k^{2/q}$,
$\beta_k=t_k^{2/q}/T_N$, and
$\mu_N=\sum_{k\le N}\beta_k\delta_{s_ke^{i\varphi_k}}$.
Set $h=\lceil q/2\rceil$ and $D_q=(h!)^{q/(2h)}$. For each fixed
$j\ge1$ the phases $e^{ij\varphi_k}$ remain independent and uniform.
Hence, by \eqref{eq:steinhaus}, Tonelli, and Minkowski in
$\ell^{q/2}$,
\begin{align*}
 \E\sum_{j\ge1}|\widehat\mu_N(j)|^q
 &\le D_q\sum_{j\ge1}
          \left(\sum_{k\le N}\beta_k^2s_k^{2j}\right)^{q/2}\\
 &\le D_q\left[
    \sum_{k\le N}\beta_k^2
       \left(\frac{s_k^q}{1-s_k^q}\right)^{2/q}\right]^{q/2}\\
 &\le D_q\left(\sum_{k\le N}\beta_k^2t_k^{-2/q}\right)^{q/2}
   =D_qT_N^{-q/2}.
\end{align*}
Choose a subsequence on which the right sides are summable. The corresponding
nonnegative moment sums are then summable almost surely, hence tend to zero.
Theorem \ref{thm:moment} rules out a nonzero vanishing witness.
\end{proof}
For one $1<p<2$, take $s_k=1-(k+1)^{-a}$ with $1<a\le q/2$:
the set is Blaschke but almost surely a uniqueness set.
For a simultaneous example take
\[
 s_k=1-\frac1{(k+2)\log^2(k+2)}.
\]
Apply the proposition to a countable sequence of exponents increasing to two and
use the inclusions of the coefficient spaces. Almost surely this gives a
Blaschke uniqueness set for every $1\le p<2$.
No converse to the random-angle criterion is claimed.

\section{Non-Blaschke exact zero sets}\label{sec:nonblaschke}
\begin{corollary}[A universal exact product for $p>2$]\label{thm:factorialproduct}
The set
\[
 V_* =\bigcup_{m\ge1}\{z\in\D:z^{m!}=e^{-1/m}\}
\]
is not Blaschke and is a zero set for $\ell_A^p$ exactly when $p>2$.
The product
\[
 F_*(z)=\prod_{m\ge1}
 \left(1-\frac{z^{m!}+z^{2m!}+\cdots+z^{mm!}}
 {e^{-1/m}+e^{-2/m}+\cdots+e^{-1}}\right)
\]
belongs to every $\ell_A^p$, $p>2$, and has exactly the points of $V_*$
as simple zeros.
\end{corollary}
\begin{proof}
Apply Theorem~\ref{thm:exactrings} with $d_m=m$. The recurrence gives
$n_m=m!$, and the criterion is $\sum_m m^{1-p}<\infty$, equivalently
$p>2$. The Blaschke criterion is $\sum_m1/m<\infty$, which fails.
The endpoint $p=1$, which is not covered by that theorem, is excluded
by $\ell_A^1\subset H^2$.
\end{proof}

\begin{remark}[A variant]\label{thm:nonblaschke}
The same conclusion holds for the alternative polygons
\[
 n_m=\frac{(m+1)!}{2},\qquad \rho_m=\frac{m}{m+1},\qquad
 V=\bigcup_{m\ge1}\{z\in\D:z^{n_m}=\rho_m\}.
\]
Indeed, with $\varepsilon_m=\log((m+1)/m)\asymp1/m$, one has
$D_m=\lceil1/\varepsilon_m\rceil=m+1$ and
$n_{m+1}=(D_m+1)n_m$. Proposition~\ref{prop:light} gives
\[
 F(z)=\prod_{m\ge1}
 \left(1-\frac{\sum_{k=1}^{m+1}z^{kn_m}}
                  {\sum_{k=1}^{m+1}(m/(m+1))^k}\right),
\]
which belongs to every $\ell_A^p$, $p>2$, and has exactly $V$ as its
simple zero set. The polygon-mass comparison in
Section~\ref{sec:obstructions} and
$\sum_{m\le M}\varepsilon_m=\log(M+1)\to\infty$ show that $V$ is not
Blaschke. Hence no nonzero function in $\ell_A^p$, $1\le p\le2$, vanishes
on $V$.
\end{remark}

\section{Dependence on the angular distribution}\label{sec:moduli}
This section shows that zero-set behavior is not determined by the moduli
alone: changing the angular placement can turn an exact zero sequence into
a uniqueness set.

\begin{proposition}[Exact zero sequences on a radius]\label{prop:radius}
Let $0<s_k<1$, with repetitions allowed, and $\sum_k(1-s_k)<\infty$.
Then the multiset of positive real points $(s_k)$ has an exact witness in
$\ell_A^1$, with precisely the specified multiplicities.
\end{proposition}
\begin{proof}
Let $B(z)=\prod_k b_{s_k}(z)$, where $b_s(z)=(s-z)/(1-sz)$, be the
Blaschke product with the prescribed zeros and multiplicities, and put
$F(z)=(1-z)^2B(z)$. Since $|b_s(z)|\le1$ in $\D$,
\[
 \frac{1-z}{1-sz}=\frac{1+b_s(z)}{1+s},\qquad
 \left|\frac{1-z}{1-sz}\right|\le\frac{2}{1+s}.
\]
Differentiating finite Blaschke products and passing to the locally
uniform limit gives
\[
 |(1-z)^2B'(z)|
 \le\sum_k(1-s_k^2)\left|\frac{1-z}{1-s_kz}\right|^2
 \le4\sum_k\frac{1-s_k}{1+s_k}<\infty.
\]
Hence $F'(z)=-2(1-z)B(z)+(1-z)^2B'(z)$ belongs to $H^\infty\subset H^2$.
Writing $F(z)=\sum_{n\ge0}a_nz^n$, Parseval and Cauchy--Schwarz give
\[
 \sum_{n\ge1}|a_n|
 \le\left(\sum_{n\ge1}\frac1{n^2}\right)^{1/2}
       \left(\sum_{n\ge1}n^2|a_n|^2\right)^{1/2}
 =\frac{\pi}{\sqrt6}\,\|F'\|_{H^2}<\infty.
\]
Thus $F\in\ell_A^1$ and $F(0)=B(0)>0$. Since $(1-z)^2$ has no zeros
in $\D$, $F$ has exactly the prescribed zeros, with exactly their
prescribed multiplicities.
\end{proof}
\begin{corollary}[Identical moduli and different zero-set behavior]\label{thm:angular}
There are two sequences in $\D$ with identical termwise moduli such that
one is an exact zero sequence for $\ell_A^1$, whereas the other is a
uniqueness set for every $\ell_A^p$, $1\le p<2$. Repetitions in the first
sequence are counted as prescribed multiplicities.
\end{corollary}
\begin{proof}
Let $(w_k)$ enumerate the polygonal uniqueness set of
Corollary~\ref{thm:dyadic}. The sequences $(|w_k|)$ and $(w_k)$, in that
order, have the required properties. The first has finite Blaschke sum,
so Proposition~\ref{prop:radius} realizes it exactly in $\ell_A^1$,
including its repetitions. The second is a uniqueness set for every
$1\le p<2$ by Corollary~\ref{thm:dyadic}.
\end{proof}

The normalized counting measure of the full polygon has sparse moments
because it is the radial image of a monomial Clark measure. Putting all
its vertices on one radius destroys that cancellation.

\section{Non-exact zero sets for \texorpdfstring{$1<p<2$}{1 < p < 2}}\label{sec:localnonexact}
We construct a Blaschke set contained in the zero set of a function in
$\bigcap_{p>1}\ell_A^p$ that is not an exact zero set for any
$1<p<2$. The proof uses quadratic characters to select polygon vertices
and shows that no function in these coefficient spaces can have precisely
the selected set as its zero set.

\subsection{Quadratic characters and Gauss sums}\label{subsec:gauss}
Let $N$ be an odd prime. A nonzero residue class $k$ modulo $N$ is a
\emph{quadratic residue} if $k\equiv a^2\pmod N$ for some integer $a$;
otherwise it is a \emph{quadratic non-residue}. The quadratic character,
or Legendre symbol, is
\[
 \chi(k)=\left(\frac{k}{N}\right)=
 \begin{cases}1,&k\text{ is a nonzero square modulo }N,\\
 -1,&k\text{ is a non-residue modulo }N,\\
 0,&N\mid k.
 \end{cases}
\]
It is regarded as a periodic function on the integers. The squaring map is two-to-one from the nonzero classes modulo $N$
onto the nonzero squares, so there
are $(N-1)/2$ residues and as many non-residues. Multiplicativity follows
from the elementary description of the multiplicative group modulo $N$,
or from Euler's criterion; see \cite{IR}. In particular $\sum_k\chi(k)=0$.
For
\[
 g_N=\sum_{a\bmod N}\chi(a)e^{2\pi ia/N}
\]
a change of variable gives
\[
 \sum_{a\bmod N}\chi(a)e^{2\pi iba/N}=\chi(b)g_N.
\]
For $N\nmid b$ this uses $\chi(b)^2=1$; at $b=0$ both sides vanish.
Parseval applied to the vector $\chi(a)$ now gives $|g_N|=\sqrt N$.
Also $\overline{g_N}=\chi(-1)g_N$, by replacing $a$ with $-a$;
hence $g_N^2=\chi(-1)N$. We will use this last identity to check the
Fourier diagonalization of the character matrix.

\subsection{Coefficient estimates for selected Clark atoms}
The construction in this section occupies the indices where $\chi(k)\ne-1$,
including the zero class; its omitted indices have $\chi(k)=-1$.

\begin{lemma}[The character matrix]\label{lem:character}
Let $n=n'N$ with $N$ an odd prime, $\omega=e^{2\pi i/n}$, and
$Y_k=\sum_{l<n}c_l\omega^{lk}$. If $Y_k=0$ whenever $\chi(k)\ne-1$, then
\[
 c_l=-\sum_{l'<n}\mathcal K(l'-l)c_{l'},\qquad
 \mathcal K(d)=\frac1n\sum_{k<n}\chi(k)\omega^{dk}.
\]
Moreover $\mathcal K(d)=0$ unless $n'\mid d$, and
$\mathcal K(n'b)=\chi(b)g_N/N$. Therefore
\begin{equation}\label{eq:charactermax}
 |c_l|\le N^{-1/2}\sum_{\substack{l'<n\\l'\equiv l\ ({\rm mod}\ n')}}|c_{l'}|.
\end{equation}
\end{lemma}
\begin{proof}
Put $G(z)=\sum_{l<n}c_lz^l\in K_{z^n}$. Its Clark transform has
values $Y_k$ at $\omega^k$. The vanishing assumption says
$\chi(k)Y_k=-Y_k$, or $\chi(U_1)G=-G$. Taking coefficients and
using~\eqref{eq:clarkmatrix} gives
$c_l=-\sum_{l'<n}\mathcal K(l'-l)c_{l'}$.

For the entries of this matrix, apply the composition identity
\eqref{eq:clarkcomposition}. On the atoms of
$\sigma_{e^{2\pi ia/N}}^{(n')}$, the value of the character is
$\chi(a)$. Hence
\[
 \mathcal K(d)=\frac1N\sum_{a=0}^{N-1}\chi(a)
       \int_\T\zeta^d\,d\sigma_{e^{2\pi ia/N}}^{(n')}(\zeta).
\]
The moments in~\eqref{eq:clarkmeasure} vanish unless $n'\mid d$.
For $d=n'b$, the last display becomes
\[
 \mathcal K(n'b)=\frac1N\sum_{a=0}^{N-1}\chi(a)e^{2\pi iab/N}
       =\frac{\chi(b)g_N}{N}.
\]
Thus $|\mathcal K(d)|\le N^{-1/2}$, and only coefficients with the same residue
class modulo $n'$ interact. Taking absolute values gives
\eqref{eq:charactermax}. This derives the block structure from
Clark-measure composition and the size estimate from the Gauss sum.
\end{proof}

\begin{lemma}[A coefficient norm estimate]\label{lem:weakspread}
For $1\le p\le2$ put
\[
 \gamma(p)=\left(\frac1p-\frac12\right)\left(1-\frac p2\right)
           =\frac{(2-p)^2}{4p}.
\]
Under the hypotheses of Lemma \ref{lem:character},
$\|c\|_2\le N^{-\gamma(p)}\|c\|_p$.
\end{lemma}
\begin{proof}
In a residue class $S$ of size $N$, put $x_S=\|(c_l)_{l\in S}\|_p$.
Equation \eqref{eq:charactermax} and H\"older give
$\max_{l\in S}|c_l|\le N^{1/2-1/p}x_S$. Hence
\[
 \sum_{l\in S}|c_l|^2
 \le N^{(1/2-1/p)(2-p)}x_S^2=N^{-2\gamma(p)}x_S^2.
\]
Sum over $S$ and use $\|(x_S)\|_2\le\|(x_S)\|_p$.
\end{proof}

The following refinement is needed only for the stronger quantitative
bound in Remark~\ref{rem:forced}; the main non-exactness proof uses
Lemma~\ref{lem:weakspread}.

\begin{lemma}[An improved estimate from a compressed projection]\label{lem:strongspread}
Let $N\ge257$ be prime, $1\le p<2$, and
\[
 \eta(p)=\frac12\left(\frac1p-\frac12\right)=\frac{2-p}{4p},
 \qquad C_p=4\cdot16^{\,1/p-1/2}\le16.
\]
If the $N$-point Fourier transform of $u$ is supported on the
quadratic non-residues, then
\begin{equation}\label{eq:strongspread}
 \|u\|_2\le C_pN^{-\eta(p)}\|u\|_p.
\end{equation}
For $n=n'N$ and the hypothesis of Lemma \ref{lem:character}, the same estimate
holds with $c$ in place of $u$.
\end{lemma}
\begin{proof}
Use the positive-exponent convention of~\eqref{eq:dft}:
\[
 (\mathcal F u)_k=\frac1{\sqrt N}\sum_{t=0}^{N-1}u_te^{2\pi itk/N}.
\]
Let $T$ be the allowed support and let $P=\mathcal F^*1_T\mathcal F$. Here $1_T$ means
multiplication by the indicator of $T$, so $P$ is the orthogonal projection
onto the vectors whose Fourier transform is supported on $T$. Then $u=Pu$.
Its diagonal entries are $|T|/N=(N-1)/(2N)$, and the Gauss sum gives
\[
 |P_{kl}|\le\frac{\sqrt N+1}{2N}\le N^{-1/2}\quad(k\ne l).
\]
Indeed $1_T=(1-\chi-1_{\{0\}})/2$.
Let $J=\lfloor\sqrt N/8\rfloor$, and let $I$ consist of the $J$ largest
coordinates of $u$. Let $P_{II}$ be the submatrix with both indices in $I$, acting on the
coordinates in $I$. Since it is Hermitian, its operator norm is at most
its largest absolute row sum. The entry bounds give
$\|P_{II}\|_{2\to2}<2/3$.
Writing $u=v+w$, with $v=1_Iu$, yields
$\|v\|_2\le(2/3)\|v\|_2+\|w\|_2$, so $\|u\|_2\le4\|w\|_2$.
Each coordinate of $w$ has modulus at most $J^{-1/p}\|u\|_p$; consequently
\[
 \|w\|_2^2\le J^{-(2-p)/p}\|u\|_p^2.
\]
Since $J\ge\sqrt N/16$, this is \eqref{eq:strongspread}.

In the composite case, fix $0\le l_0<n'$ and put
$x_t=c_{l_0+n't}$, $0\le t<N$. Lemma~\ref{lem:character} gives
$x=-Mx$, where $M_{t,t'}=\chi(t'-t)g_N/N$. With
$\widehat x(k)=\sum_{t=0}^{N-1}x_te^{2\pi itk/N}$, the Gauss-sum
identity and $g_N^2=\chi(-1)N$ give
\begin{align*}
 \widehat{Mx}(k)
 &=\frac{g_N}{N}\left(\sum_{d=0}^{N-1}\chi(d)e^{-2\pi idk/N}\right)
     \widehat x(k)\\
 &=\frac{g_N^2}{N}\chi(-k)\widehat x(k)
 =\chi(k)\widehat x(k).
\end{align*}
Thus $(1+\chi(k))\widehat x(k)=0$: the transform of each class vector
is supported on the non-residues. Apply~\eqref{eq:strongspread} to
each class and sum squares as in the preceding lemma.

Here $P$ is the spectral projection of the Clark unitary onto the allowed
atoms. The compression $P_{II}$ tests whether a vector in that spectral
subspace could concentrate on very few monomials. Its norm being strictly
less than one prevents that concentration.
\end{proof}
The gain $\eta(p)$ is stronger than $\gamma(p)$ for $1<p<2$.
The estimate with $\gamma(p)$ has constant one and gives a shorter local
proof. The estimate with $\eta(p)$ strengthens the quantitative
lower bound for the norm of an exact witness, as explained in
Remark~\ref{rem:forced}.

\subsection{The construction}
Choose increasing odd primes $N_m$ with $\log N_m\ge4^m$, and put
\[
 n_0=1,\quad n_m=N_1\cdots N_m,\quad
 \varepsilon_m=2^{-m},\quad \rho_m=e^{-\varepsilon_m},\quad
 r_m=\rho_m^{1/n_m},\quad\omega_m=e^{2\pi i/n_m}.
\]
For definiteness each prime can be chosen as the least one greater than
both its predecessor and $e^{4^m}$.
With $\chi_m$ the character modulo $N_m$, define
\[
 \mathcal Q_m=\{k\in\mathbb Z_{n_m}:\chi_m(k)\ne-1\},\qquad
 W_m=\{r_m\omega_m^k:k\in \mathcal Q_m\},\qquad W=\bigcup_mW_m.
\]
The full polygon
$R_m=\{z:z^{n_m}=\rho_m\}$ is the preimage under
$z\mapsto z^{n_{m-1}}$ of the full polygon
$\{u:u^{N_m}=\rho_m\}$. Its subset $W_m$ is the preimage of the
vertices indexed by the nonzero quadratic residues together with zero.
The radii strictly increase to one.

The full polygon is the radial image of $\sigma_1^{(n_m)}$.
Under \eqref{eq:clarkcomposition}, the selected atoms are exactly those
blocks $\sigma_\beta^{(n_{m-1})}$ for which
$\beta=e^{2\pi ia/N_m}$ has $\chi_m(a)\ne-1$.

\begin{theorem}\label{thm:localnonexact}
The set $W$ above is Blaschke and is a zero set for every $\ell_A^p$, $p>1$.
It is not an exact zero set for any $1<p<2$, as a set, with multiplicities ignored. It is exact for $p\ge2$.
There is a single function $F\in\bigcap_{p>1}\ell_A^p$ with only
simple zeros such that $Z(F)=\bigcup_mR_m\supset W$.
\end{theorem}
\begin{proof}[Proof of Theorem~\ref{thm:localnonexact}: existence]
Apply Proposition~\ref{prop:light} with $D_m=d_m=2^m$ and
$\rho_m=e^{-1/d_m}$; its separation condition holds because
$n_{m+1}=N_{m+1}n_m\ge(d_m+1)n_m$.
The resulting product $F$ has exactly the full polygons as its simple zeros.
As in the proof of Theorem~\ref{thm:exactrings}, the estimate
$\sum_{k=1}^{d_m}e^{-k/d_m}\ge d_m/e$ gives
\begin{equation}\label{eq:localwitnessnorm}
 \|F\|_p^p\le
 \exp\left(e^p\sum_m2^{-m(p-1)}\right)
 =\exp\left(\frac{e^p}{2^{p-1}-1}\right).
\end{equation}
The moment theorem therefore gives the explicit lower bound
\[
 \delta_q(W)\ge
 \left[\exp\left(\frac{e^p}{2^{p-1}-1}\right)-1\right]^{-1/p}.
\]
Also $\sum_{w\in W}(1-|w|)\le\sum_m\varepsilon_m=1$. The Blaschke
product of $W$ gives exact realization for every $p\ge2$.
\end{proof}
Non-exactness is proved in Section~\ref{subsec:local-completion}, using
the local estimate below. We isolate that estimate because the endpoint
construction uses the same mechanism.

\subsection{A one-level test for the absence of extra zeros}
\begin{proposition}[A local estimate in zero-free regions]\label{prop:boxes}
Let $n=n'N$, where $N$ is an odd prime. Put
\begin{gather*}
 \rho\in(1/3,1),\qquad \tau=\frac{1-\rho}{2},\qquad
 s=\frac{1+\rho}{2},\qquad \Delta=D(\rho,\tau),\\
 \omega=e^{2\pi i/n},\qquad
 \Gamma=\{\alpha\in\T:|\alpha-1|\le\tau/(2\rho)\}.
\end{gather*}
Let $f$ be analytic in $\D$, with $f(0)=1$, and suppose
\begin{equation}\label{eq:boxupper}
 \frac1n\sum_{z^n=x}|f(z)|^2\le B^2
       \quad(0<|x|<s)
\end{equation}
for some $B\ge1$. Set
$\mathcal Q=\{k\in\mathbb Z_n:\chi(k)\ne-1\}$ and
$E=\mathbb Z_n\setminus\mathcal Q$. For every $k\in\mathbb Z_n$,
define $\Phi_k(x)=f(x^{1/n}\omega^k)$ on $\Delta$, using the
principal root. Assume that $\Phi_e$ is zero-free on $\Delta$ for
every $e\in E$. Thus $f$ has no
zeros in the regions $\{x^{1/n}\omega^e:x\in\Delta\}$ for $e\in E$;
we call these regions the \emph{omitted boxes}.
Finally suppose
\[
 \left(\frac1{|E|}\sum_{e\in E}|\Phi_e(\rho)|^2\right)^{1/2}
        \le\sqrt3\,A N^{-\lambda}B
\]
with $A\ge1$, $\lambda>0$. Then
\begin{equation}\label{eq:boxconclusion}
 \log B\ge\frac{\delta_\Gamma\lambda}{15}\log N
                  -\log(\sqrt{48}\,A^{1/5}),\qquad
 \delta_\Gamma=\frac2\pi\arcsin\frac{\tau}{4\rho}.
\end{equation}
\end{proposition}
No vanishing on the occupied vertices is assumed in this proposition.
In its applications, such vanishing is used separately to obtain the
central sample bound from the coefficient estimates above.
\begin{proof}
There are $|E|=n(N-1)/(2N)\ge n/3$ omitted indices and $|\mathcal Q|\ge n/2$
occupied ones. All boxes are well defined and disjoint: $\Delta$ is in
the right half-plane and taking the principal $n$-th root changes the
argument by less than $\pi/(2n)$.

\paragraph{Local maxima.}
Put $M_k=\sup_{|x-\rho|<3\tau/4}|\Phi_k(x)|$.
For every point in this smaller disk, the disk of radius $\tau/4$ about it
lies in $\Delta$. The area sub-mean inequality and \eqref{eq:boxupper} give
\[
 M_k^2\le\frac{16}{\pi\tau^2}\int_\Delta|\Phi_k(x)|^2\,dA(x),
 \qquad \sum_k M_k^2\le16nB^2.
\]
Consequently
\[
 \frac1{|E|}\sum_{e\in E}\log M_e\le\log(\sqrt{48}B),\qquad
 \frac1{|\mathcal Q|}\sum_{k\in \mathcal Q}\log M_k\le\log(\sqrt{32}B).
\]
Concavity of the logarithm also gives
\[
 \frac1{|E|}\sum_{e\in E}\log|\Phi_e(\rho)|
       \le\log(\sqrt3 AB)-\lambda\log N.
\]

\paragraph{Propagation through a zero-free box.}
For $e\in E$, $\log M_e-\log|\Phi_e|$ is nonnegative and harmonic on
$D(\rho,3\tau/4)$. Harnack, at distance at most $\tau/2$, gives
\[
 \log|\Phi_e(x)|
 \le\frac45\log M_e+\frac15\log|\Phi_e(\rho)|
       \quad(|x-\rho|\le\tau/2).
\]
Average over $E$, use $\sqrt3\le\sqrt{48}$, and combine with the bound
over $\mathcal Q$. The negative term has weight $|E|/n\ge1/3$, so
\begin{equation}\label{eq:localboxlog}
 \frac1n\sum_{k<n}\log|\Phi_k(x)|
 \le\log B+\log(\sqrt{48}A^{1/5})
                    -\frac{\lambda}{15}\log N
       \quad(|x-\rho|\le\tau/2).
\end{equation}

\paragraph{Averaging around the circle.}

Let $L=L_f$ be the logarithmic average from
Section~\ref{subsec:logarithmic-averages}, with the present $n$ and $\rho$.
By~\eqref{eq:clark-log-average}, $\int_\T L\dm\ge0$.
The bound~\eqref{eq:boxupper} gives $L\le\log B$ everywhere, while
\eqref{eq:localboxlog} applies on $\Gamma$, whose normalized length is
$\delta_\Gamma$. Integrating these bounds yields
\[
 0\le\log B+\delta_\Gamma\log(\sqrt{48}A^{1/5})
                  -\frac{\delta_\Gamma\lambda}{15}\log N.
\]
Since $\delta_\Gamma\le1$ and $\log(\sqrt{48}A^{1/5})\ge0$, this implies
\eqref{eq:boxconclusion}. No prescribed multiplicity was used.
\end{proof}

\subsection{Application to the infinite set}\label{subsec:local-completion}
\begin{proof}[Completion of the proof of Theorem~\ref{thm:localnonexact}]
Fix $1<p<2$ and suppose $f$ is an exact witness, normalized by $f(0)=1$.
Put $C=\|f\|_p$, and at level $m$ abbreviate
$n=n_m$, $N=N_m$, $\rho=\rho_m$, $\varepsilon=2^{-m}$, and
$\mathcal Q=\mathcal Q_m$.
The fold estimate and $\ell^p\subset\ell^2$ give
\eqref{eq:boxupper} with
\[
 B=\kappa_p(s)C,\qquad s=(1+\rho)/2.
\]
At the central polygon, Lemma \ref{lem:weakspread} gives the sample bound in
Proposition \ref{prop:boxes} with $A=1$ and $\lambda=\gamma(p)$.

We verify that each omitted box contains no prescribed zero. If $w$ is on level $j<m$, then
\[
 |w|^n=e^{-\varepsilon_jn/n_j}\le e^{-\varepsilon N}\le e^{-2}<\rho-\tau.
\]
Indeed $N\ge2^{m+1}$ and $\rho-\tau=(3\rho-1)/2>0.4$.
If $j>m$, then
\[
 |w|^n=e^{-\varepsilon_jn/n_j}
 \ge1-\varepsilon_j/N_{m+1}>1-\varepsilon/4\ge s,
\]
because $(1-\rho)/2\ge\varepsilon/4$.
On the present polygon the only prescribed zeros have indices in $\mathcal Q$.
Thus no prescribed zero lies in an omitted box; exactness makes every
$\Phi_e$ zero-free there. This is the only use of exactness.

Finally
\[
 \delta_\Gamma\ge\frac{\varepsilon}{8\pi},\qquad
 \kappa_p(s)\le(1-s)^{-1/q}\le(4/\varepsilon)^{1/q}.
\]
Equation \eqref{eq:boxconclusion} yields
\begin{equation}\label{eq:localcontradiction}
 \log C\ge\frac{\gamma(p)}{120\pi}\,2^m
                  -\frac{m+2}{q}\log2-\log\sqrt{48}\longrightarrow\infty.
\end{equation}
This contradiction proves Theorem~\ref{thm:localnonexact}.
\end{proof}

\begin{remark}[Stronger constants and forced extra zeros]\label{rem:forced}
For $m\ge2$, one has $N_m\ge257$. Lemma~\ref{lem:strongspread}
then permits replacing $\gamma(p)$ by $\eta(p)$ in
\eqref{eq:localcontradiction}, subtracting the additional constant
$\frac15\log C_p$.
More importantly, the argument is local. Every normalized vanishing witness
has an unwanted zero in at least one omitted box at every sufficiently large
level. The general conclusions about unwanted zeros are those of
Lemma~\ref{lem:division}.
\end{remark}

\section{Non-exact zero sets in the analytic Wiener algebra}\label{sec:endpoint}
The full-polygon construction of Section~\ref{sec:localnonexact}
accumulates at every point of $\T$. It cannot be a zero set for
$\ell_A^1$, whose functions are continuous on the closed disk: such a
function would vanish on the whole boundary and hence identically.
We therefore localize the angular argument near a single boundary point.

\subsection{Character sums on intervals and localization}
We first record the arithmetic estimate needed for localization. With
$N$ an odd prime and $\chi$ its quadratic character as defined in
Section~\ref{subsec:gauss}, one has the standard character-sum bound
\begin{equation}\label{eq:pv}
 \left|\sum_{u<k\le v}\chi(k)\right|\le\sqrt N\log N.
\end{equation}
Here one can discard full periods and assume that the interval has length
at most $N$. Fourier inversion and the Gauss-sum identity bound its sum by
\[
 \frac{\sqrt N}{N}\sum_{b=1}^{N-1}
   \min\left(v-u,\frac1{|\sin(\pi b/N)|}\right)
 \le\sqrt N\sum_{b=1}^{(N-1)/2}\frac1b
 \le\sqrt N\log N.
\]
For the last inequality integrate $1/x$ over the intervals
$[b-\tfrac12,b+\tfrac12]$: convexity gives
$1/b\le\int_{b-1/2}^{b+1/2}dx/x$. For more general character-sum
estimates see~\cite[Chapter~23]{Da}.

\subsubsection{Fej\'er's identity}
For each integer $K\ge1$,
\[
 \sum_{|s|<K}\left(1-\frac{|s|}{K}\right)e^{ist}
   =\frac1K\left|\sum_{u=0}^{K-1}e^{iut}\right|^2\ge0.
\]
Indeed, exactly $K-|s|$ pairs $(u,u')$ have $u-u'=s$.
The nonnegativity will be used to construct a finitely supported
weight whose Fourier coefficients have $\ell^1$ norm one.

\subsubsection{A local estimate for the omitted indices}
\begin{lemma}\label{lem:localcharacter}
Let $N$ be an odd prime and
$3\sqrt N\log N+3\le K\le N/4$, with $K$ integral.
Suppose
\[
 Y_k=\sum_{l\bmod N}C_le^{2\pi ilk/N},\qquad \sum_l|C_l|\le C_0,
\]
and $|Y_k|\le\varepsilon_{\mathrm{err}} C_0$ whenever $|k-k_0|<K$ and
$\chi(k)\ne-1$, where $\varepsilon_{\mathrm{err}}\sqrt{2K}\le1$.
Let $E$ consist of the integers $-K/2\le k-k_0<K/2$ with $\chi(k)=-1$.
Then
\[
 |E|\ge K/3,\qquad
 \frac1{|E|}\sum_{k\in E}|Y_k|^2\le\frac{48\sqrt N}{K}C_0^2.
\]
\end{lemma}
\begin{proof}
The middle interval contains $K$ integers, at most one a multiple of $N$.
The character-sum bound \eqref{eq:pv} gives
$|E|\ge(K-1-\sqrt N\log N)/2\ge K/3$.

Periodize the tent $w(k)=\max(0,1-|k-k_0|/K)$ modulo $N$. Its Fourier
coefficients satisfy
\[
 \widehat w(l)=\frac{e^{-2\pi ilk_0/N}}N
       \sum_{|s|<K}(1-|s|/K)e^{-2\pi ils/N}.
\]
Fej\'er's identity shows that
$\widehat w(l)e^{2\pi ilk_0/N}\ge0$.
Inversion at $k_0$ therefore gives $\sum_l|\widehat w(l)|=1$.
The coefficients $d$ of $Z=Yw$ consequently satisfy $\|d\|_1\le C_0$.

Write $Z=Z'+Z''$, where $Z'$ is supported on the non-residues and $Z''$
on the other positions. The latter has fewer than $2K$ nonzero entries,
each at most $\varepsilon_{\mathrm{err}} C_0$. Parseval and finite norm comparison give
\[
 \|d''\|_1\le\sqrt N\|d''\|_2
       =\left(\sum_k|Z''_k|^2\right)^{1/2}\le\sqrt{2K}\varepsilon_{\mathrm{err}} C_0\le C_0.
\]
Thus $\|d'\|_1\le2C_0$. Apply Lemma~\ref{lem:weakspread} with
$p=1$ and $n=N$ to $d'$. Since $\gamma(1)=1/4$, Parseval gives
\[
 \sum_k|Z'_k|^2=N\|d'\|_2^2
      \le\sqrt N\,\|d'\|_1^2\le4\sqrt N\,C_0^2.
\]
For $k\in E$, $w(k)\ge1/2$ and $Z'_k=Y_kw(k)$, so
$\sum_{k\in E}|Y_k|^2\le16\sqrt N C_0^2$.
Divide by $|E|\ge K/3$.

Multiplication by $w$ on the Clark atoms is $w(U_1)$.
The identity $\sum|\widehat w|=1$ says that this circulant matrix is
contractive for coefficient $\ell^1$ norms; this is the localization step.
\end{proof}

\subsection{A lacunary series and the non-exact subset}
\begin{theorem}\label{thm:endpoint}
Let $N_1<N_2<\cdots$ be primes such that
\[
 \log N_1\ge100,\qquad
 \log N_j\ge10^{10}jN_{j-1}\quad(j\ge2).
\]
Define
\begin{equation}\label{eq:endpointseries}
 f(z)=1-\frac32\sum_{j\ge1}(-1/2)^{j-1}z^{N_j}
      =1-\frac32z^{N_1}+\frac34z^{N_2}-\frac38z^{N_3}+\cdots.
\end{equation}
Then $f\in\ell_A^1$, $\|f\|_1=4$, $f(0)=1$, and $f(1)=0$.
The subset $W^{(1)}\subset Z(f)$ defined below accumulates only at $1$ and is
not an exact zero set for $\ell_A^1$. It is nevertheless an exact zero set
for $H^\infty$ and for the disk algebra.
\end{theorem}

The proof occupies Sections~\ref{subsec:endpoint-locate}--\ref{subsec:endpoint-propagation}.

\subsection{Locating the zeros}\label{subsec:endpoint-locate}
\noindent\emph{Proof of Theorem~\ref{thm:endpoint}.}
Write $c_j=(3/2)(-1/2)^{j-1}$. Then $\sum|c_j|=3$ and $\sum c_j=1$.
The prime-growth hypotheses imply, for $m\ge2$,
\begin{equation}\label{eq:endpointgrowth}
 N_{m+1}\ge5N_m^2,\qquad
 N_m^{0.1}\ge3000N_{m-1},\qquad \log N_m\ge100,
\end{equation}
and $N_{m-1}\ge2^mN_{m-2}$ for $m\ge3$.
For example
$0.1\log N_m\ge10^9mN_{m-1}\ge\log3000+\log N_{m-1}$;
the other inequalities follow even more directly.

Fix $m\ge2$, and abbreviate $n=N_m$, $M=N_{m-1}$. Here the local
polygon size $n$ is a single prime, rather than the product of primes
used in Section~\ref{sec:localnonexact}. Set
\[
 x(z)=c_m^{-1}\left(1-\sum_{j<m}c_jz^{N_j}\right),\qquad
 H(z)=c_m^{-1}\sum_{j>m}c_jz^{N_j}.
\]
Thus $f(z)=c_m(x(z)-z^n-H(z))$ and $x(1)=2/3$.
Since $|c_j/c_m|=2^{m-j}$,
\[
 |x'(z)|\le\sum_{j<m}2^{m-j}N_j
       \le2M+2^mN_{m-2}\le4M\quad(|z|\le1),
\]
where the earlier-index term is omitted when $m=2$.
Consequently
\begin{equation}\label{eq:xnear}
 |x(z)-2/3|\le1/24
 \quad\text{if }|z|\le1,\quad |z-1|\le(96M)^{-1}.
\end{equation}
For $j>m$, $N_j/n\ge5n$ and $\sum_{j>m}2^{m-j}=1$, so
\begin{equation}\label{eq:Hsmall}
 |H(z)|\le(4/5)^{5n}\le e^{-n}\quad\text{if }|z|^n\le4/5.
\end{equation}

\begin{lemma}\label{lem:endpointzeros}
Let $p_k=(2/3)^{1/n}e^{2\pi ik/n}$ and
$k_m=\lfloor n/(1000M)\rfloor$. For each integer $|k|\le k_m$ there
is exactly one zero $w_k$ of $f$ in $|z-p_k|<1/(8n)$.
It is simple and satisfies
\[
 |w_k^n-2/3|\le1/24+e^{-n},\qquad
 w_k=(w_k^n)^{1/n}e^{2\pi ik/n},
\]
with the principal root on the right.
\end{lemma}
\begin{proof}
For $|z-p_k|\le1/(8n)$,
\[
 |z-1|\le\frac{0.6}{n}+\frac{2\pi}{1000M}\le\frac1{96M},
 \qquad
 |z|^n\le\frac23e^{1/7}<\frac45.
\]
In deriving the second bound use
$|z|\le(2/3)^{1/n}(1+1/(7n))$; this also puts the whole disk inside $\D$.
Thus \eqref{eq:xnear}--\eqref{eq:Hsmall} apply.

On the boundary write $z=p_k(1+v)$. Then
$1/(8n)\le|v|\le1/(7n)$. Put $u=n\log(1+v)$.
The convergent logarithm series gives
$|u-nv|\le n|v|^2\le0.001$, so $0.124\le|u|\le0.144$.
The Taylor expansion gives
$|e^u-1|\ge|u|-\frac12|u|^2e^{|u|}>0.11$. Hence
\[
 |z^n-2/3|=\frac23|e^u-1|>0.073
       >1/24+e^{-n}\ge|x(z)-2/3|+|H(z)|.
\]
Rouch\'e compares $z^n-x(z)+H(z)=-f(z)/c_m$ with $z^n-2/3$.
The latter has exactly its single simple zero $p_k$ in the disk, so the
same count, including multiplicity, holds for $f$.

At $w=w_k$ one has $w^n=x(w)-H(w)$, proving the first estimate.
Put $y=w^n$; then $|(3/2)y-1|<0.0626$ and
$|\log((3/2)y)|<0.07$. Thus
$|y^{1/n}e^{2\pi ik/n}-p_k|<0.071/n$.
This chosen root is within $0.2/n$ of $w$, while distinct $n$-th roots
of $y$ are at least $3/n$ apart. Since $w$ is one of them, it is the
chosen root.
\end{proof}

At level $m$ use these zeros and the character $\chi_m$ modulo $N_m$:
\begin{equation}\label{eq:endpointW}
 W_m=\{w_k:|k|\le k_m,\ \chi_m(k)\ne-1\},\qquad W^{(1)}=\bigcup_{m\ge2}W_m.
\end{equation}
Every such point is a zero of the fixed function \eqref{eq:endpointseries},
and
\begin{equation}\label{eq:endpointradii}
 0.62\le |w|^{N_m}\le0.71\quad(w\in W_m).
\end{equation}
Different levels are radially disjoint by \eqref{eq:endpointgrowth}.
Their points lie within $(96N_{m-1})^{-1}$ of $1$. Also
\[
 \sum_{w\in W^{(1)}}(1-|w|)
 \le\sum_m(2k_m+1)\frac{0.5}{N_m}<\infty,
\]
using $-\log0.62<0.5$. This proves the zero-set and Blaschke assertions
and the stated accumulation property.

\subsection{Local polygons and approximate vanishing}
Suppose $g\in\ell_A^1$ has exactly the zeros $W^{(1)}$, with $g(0)=1$.
Let $C=\|g\|_1\ge1$, so $|g|\le C$ in the closed disk.
Fix $m$, and retain $n,M,k_m$.
Put $K=\lfloor n^{0.6}\rfloor$. The growth conditions imply
\[
 3\sqrt n\log n+3\le K\le n/4,\qquad k_m\ge3K.
\]
For integers $i$ with $|i|\le\lfloor k_m/K\rfloor-1$ set
\[
 b_i=iK,\qquad y_i=w_{b_i}^n,\qquad
 q_k^{(i)}=y_i^{1/n}e^{2\pi ik/n}.
\]
The zero $w_{b_i}$ is defined even if $b_i$ is not an occupied index.
There are at least $k_m/K$ such blocks. Let
\[
 E_i=\{k:-K/2\le k-b_i<K/2,\ \chi_m(k)=-1\}.
\]
These sets are disjoint, $|E_i|\ge K/3$ by
Lemma \ref{lem:localcharacter}, and $\sum_i|E_i|\ge k_m/3$.

If $|k-b_i|<K$ and $\chi_m(k)\ne-1$, then $w_k\in W^{(1)}$, so $g(w_k)=0$.
The points $w_k^n$ and $y_i$ lie in $|y-2/3|<1/23$, where
$|(d/dy)y^{1/n}|\le1.7/n$.
Using $w^n=x(w)-H(w)$ at the zeros and $|x'|\le4M$,
\[
 |w_k^n-y_i|
 \le4M|w_k-w_{b_i}|+2e^{-n}
 \le4M\left(\frac{2\pi K}{n}+\frac1{4n}\right)+2e^{-n}
 \le\frac{29MK}{n}.
\]
It follows that
\[
 |w_k-q_k^{(i)}|\le\frac{50MK}{n^2}.
\]
Both points have modulus at most $R=0.71^{1/n}$. Since
$g(z)=\sum a_jz^j$ and $\sum|a_j|=C$,
\[
 |g'(z)|\le C\max_{j\ge1}jR^{j-1}
       \le\frac{C}{eR\log(1/R)}\le1.1Cn\quad(|z|\le R).
\]
The segment joining the points stays in this disk. Therefore
\begin{equation}\label{eq:localapprox}
 |g(q_k^{(i)})|\le55CMKn^{-1}\le\varepsilon_{\mathrm{err}} C,\qquad
 \varepsilon_{\mathrm{err}}=55Mn^{-0.4},\qquad
 \varepsilon_{\mathrm{err}}\sqrt{2K}\le78Mn^{-0.1}\le1.
\end{equation}

Folding with $\xi=y_i^{1/n}$ and $p=1$ gives
$g(q_k^{(i)})=\sum_l C_le^{2\pi ilk/n}$ with $\sum|C_l|\le C$.
Apply Lemma \ref{lem:localcharacter} using \eqref{eq:localapprox}:
\[
 \frac1{|E_i|}\sum_{k\in E_i}|g(q_k^{(i)})|^2
       \le\frac{48\sqrt n}{K}C^2\le49n^{-0.1}C^2.
\]
Thus
\begin{equation}\label{eq:endpointlogsmall}
 \frac1{|E_i|}\sum_{k\in E_i}\log|g(q_k^{(i)})|
       \le\log C+\log7-0.05\log n.
\end{equation}

\subsection{Zero-free boxes, propagation, and Jensen}\label{subsec:endpoint-propagation}
For $k\in E_i$, put
\[
 \Delta_i=D(y_i,1/6),\qquad
 \Phi_k(x)=g(x^{1/n}e^{2\pi ik/n}).
\]
The disks lie in $\{\Re x>0,\ 0.45<|x|<0.88\}$, so the principal roots
are analytic there. No prescribed zero lies in this box.
Indeed, a prescribed zero on the same level is $w_{k'}$ with
$\chi_m(k')\ne-1$, whereas $\chi_m(k)=-1$. By
Lemma~\ref{lem:endpointzeros},
$w_{k'}=(w_{k'}^n)^{1/n}e^{2\pi i k'/n}$, with the principal root.
If also $w_{k'}=x^{1/n}e^{2\pi i k/n}$ for $\Re x>0$, then
$x=w_{k'}^n$ and $e^{2\pi i k'/n}=e^{2\pi i k/n}$. Hence
$k'\equiv k\pmod n$, which is impossible because
$|k|,|k'|\le k_m<n/2$ and $k'\ne k$. At an earlier level $j<m$,
$|w|^n\le0.71^{n/N_j}<0.45$; at a later level $j>m$,
$|w|^n\ge0.62^{n/N_j}>0.88$.
Exactness therefore makes $\Phi_k$ zero-free in $\Delta_i$.

As in the propagation step of the proof of Proposition~\ref{prop:boxes},
apply Harnack to $\log C-\log|\Phi_k|$, now with radius $1/6$ and
distance at most $1/12$. This gives
\begin{equation}\label{eq:endpointHarnack}
 \log|\Phi_k(x)|\le\frac23\log C+\frac13\log|g(q_k^{(i)})|
              \quad(|x-y_i|\le1/12).
\end{equation}
Let $r=(2/3)^{1/n}$ and
$\Gamma=\{\alpha\in\T:|\alpha-1|\le1/48\}$.
Its normalized length is
$\delta_\Gamma=(2/\pi)\arcsin(1/96)\ge1/(48\pi)$.
For $\alpha\in\Gamma$,
\[
 |(2/3)\alpha-y_i|\le1/72+1/24+e^{-n}<1/12.
\]

Repeat the final averaging step of the proof of
Proposition~\ref{prop:boxes}. The indices in $\bigcup_iE_i$ are distinct,
each has Clark mass $1/n$, and~\eqref{eq:endpointHarnack} applies on
$\Gamma$; all other values are bounded by $\log C$.
Thus~\eqref{eq:clark-log-average} gives
\begin{align*}
 0&\le\log C-\frac{\delta_\Gamma}{3n}
       \sum_i\sum_{k\in E_i}\bigl(\log C-\log|g(q_k^{(i)})|\bigr)\\
  &\le\log C-\frac{\delta_\Gamma k_m}{9n}(0.05\log n-2),
\end{align*}
where the second line uses~\eqref{eq:endpointlogsmall},
$\sum_i|E_i|\ge k_m/3$, and $\log7<2$.
The floor estimate $k_m/n\ge1/(1001M)$ follows from the very large
$n/M$. Since $9\cdot48\pi\cdot1001<1.4\cdot10^6$,
\[
 \log C\ge\frac{0.05\log n-2}{1.4\cdot10^6M}\ge300m,
\]
where $\log n\ge10^{10}mM$ was used. Let $m\to\infty$, a contradiction.

This proves non-exactness, with no restriction on the multiplicities
of $g$. Since $W^{(1)}$ is Blaschke, its Blaschke product $B_{W^{(1)}}$
is an exact bounded witness. Its only possible boundary singularity is at
$1$, so $(1-z)B_{W^{(1)}}$ is continuous on the closed disk and exact in the disk
algebra. This does not imply that it has summable coefficients.
This completes the proof of Theorem~\ref{thm:endpoint}.\qed

As in Remark~\ref{rem:forced}, every normalized witness has an extra
zero in an omitted box at every sufficiently large level.

\section{Two exact zero sets whose union is a uniqueness set}\label{sec:unions}
This section concerns closure under unions, not the distinction
between a zero set and an exact zero set. Each of the two sets below has an
exact witness, but their union admits no nonzero vanishing witness at all.
At $p=2$ this cannot happen: the union of two Blaschke sets is Blaschke,
and its Blaschke product is an exact witness in $H^2$.

\begin{theorem}[Failure of union closure for $p>2$]\label{union-thm:union-failure}
For every $p>2$, there exists an exact zero set $W$ for $\lp$ such that $W\cap(-W)=\varnothing$ and $W\cup(-W)$ is a uniqueness set for $\lp$. All zeros in the two exact witnesses can be taken to be simple.
\end{theorem}

\subsection{The two estimates behind the construction}\label{subsec:union-estimates}
For $p>2$, we use a Gaussian series $F_s$ with
\[
 \frac1p<s<\frac14+\frac1{2p}.
\]
The lower endpoint makes its coefficients belong to $\ell^p$ almost surely;
the upper endpoint makes the combined zero growth exceed the allowable
bound. We prove both estimates and then the theorem. Throughout this
section, a \emph{standard complex Gaussian} has density
$\pi^{-1}e^{-|z|^2}$ with respect to area measure on $\mathbb C$;
in particular, its second absolute moment is one.

For a nonzero analytic function $h(z)=\sum_{n\ge m}a_nz^n$, with $a_m\ne0$, write
\[
 J_h(r)=\sum_{\substack{h(w)=0\\0<|w|<r}}\log\frac r{|w|}
 =\frac1{2\pi}\int_0^{2\pi}\log|h(re^{it})|\,dt-m\log r-\log|a_m|,
 \qquad 0<r<1,
\]
where zeros are counted with multiplicity. The integral is well defined even when the circle contains a zero; zeros on the circle contribute zero to the Jensen sum.

\begin{lemma}\label{union-lem:jensen-upper}
Let $p>2$ and $\beta_p=\frac12-\frac1p$. If $h\in\lp$ and $h(0)\ne0$, then
\[
 J_h(r)\le\beta_p\log\frac1{1-r}+C_h
 \qquad(0<r<1)
\]
for a finite constant $C_h$.
\end{lemma}
\begin{proof}
Write $h(z)=\sum_{n\ge0}a_nz^n$. Parseval and H\"older give
\begin{align*}
 M_2(h,r)^2
 &:=\frac{1}{2\pi}\int_0^{2\pi}|h(re^{it})|^2\,dt
 =\sum_{n\ge0}|a_n|^2r^{2n}\\
 &\le\|h\|_p^2
 \left(\sum_{n\ge0}r^{2np/(p-2)}\right)^{(p-2)/p}.
\end{align*}
Since $2p/(p-2)>1$, it follows that
\[
 M_2(h,r)\le\|h\|_p
 (1-r^{2p/(p-2)})^{-\beta_p}
 \le\|h\|_p(1-r)^{-\beta_p}.
\]
The arithmetic mean--geometric mean inequality in Jensen's formula yields
$J_h(r)\le\log M_2(h,r)-\log|h(0)|$, proving the assertion.
\end{proof}

\begin{lemma}\label{union-lem:gaussian-jensen}
Let $0<s<\frac12$, and let $\xi_n$, $n\ge0$, be independent standard
complex Gaussian random variables. Set
\[
 F_s(z)=\sum_{n\ge0}\frac{\xi_n}{(n+1)^s}z^n,
 \qquad r_k=1-e^{-k^2}.
\]
Almost surely, $F_s$ is analytic in $\D$, its zeros are simple, and
\[
 \lim_{k\to\infty}\frac{J_{F_s}(r_k)}{k^2}
 =\frac12-s.
\]
Almost surely its zero set is disjoint from its reflection through the origin.
\end{lemma}
\begin{proof}
For each $0<r<1$, the expectation of
$\sum_n|\xi_n|(n+1)^{-s}r^n$ is finite. Applying this on a countable sequence of radii increasing to $1$ gives almost-sure local absolute convergence. Also $F_s(0)=\xi_0\ne0$ almost surely.

The variance at any point of modulus $r$ is
\[
 V_s(r)=\E|F_s(re^{it})|^2
 =\sum_{n\ge0}\frac{r^{2n}}{(n+1)^{2s}}
 \asymp(1-r)^{-(1-2s)}\qquad(r\to1).
\]
For the lower bound, sum over $n\le(1-r)^{-1}$, on which $r^{2n}$ is bounded below uniformly for $r$ close to $1$. For the upper bound, use $r^{2n}\le e^{-2n(1-r)}$ and compare the resulting sum with an integral. Thus
\[
 \tfrac12\log V_s(r)=\left(\tfrac12-s\right)\log\frac1{1-r}+O(1).
\]
Define
\[
 A_s(r)=\frac{1}{2\pi}\int_0^{2\pi}
 \log\left|\frac{F_s(re^{it})}{V_s(r)^{1/2}}\right|\,dt.
\]
For each fixed $r,t$, the normalized variable in this integral is standard complex Gaussian. Since $\E|\log|\xi_0||<\infty$, Tonelli gives the uniform bound
\[
 \E|A_s(r)|\le\E|\log|\xi_0||=:C<\infty.
\]
For each $\varepsilon>0$, Markov's inequality gives
\[
 \sum_{k\ge1}\mathbb P\{|A_s(r_k)|>\varepsilon k^2\}
 \le\frac C\varepsilon\sum_{k\ge1}k^{-2}<\infty.
\]
By Borel--Cantelli, first for $\varepsilon=1/j$ and then for all $\varepsilon>0$, we have $A_s(r_k)/k^2\to0$ almost surely. Jensen's formula now gives
\[
 J_{F_s}(r_k)=\tfrac12\log V_s(r_k)+A_s(r_k)-\log|\xi_0|,
\]
which proves the limit.

To prove simplicity, condition on all coefficients except $\xi_0$. The derivative $F_s'$ is then a nonzero analytic function almost surely and has countably many zeros in $\D$. A multiple zero of $F_s$ would require $\xi_0$ to equal one of countably many specified values. Its continuous distribution makes this event have probability zero.

Finally, write $F_s(z)=E_s(z^2)+zO_s(z^2)$, where $E_s$ uses the even-indexed coefficients and $O_s$ the odd-indexed coefficients. These are independent nonzero random analytic functions. If $w\ne0$ and $F_s(w)=F_s(-w)=0$, then $E_s(w^2)=O_s(w^2)=0$. Conditional on $E_s$, its zeros are countable, and at each of them $O_s$ has a nondegenerate Gaussian value. Thus $E_s$ and $O_s$ have no common zero almost surely. The case $w=0$ has already been excluded.
\end{proof}

\subsection{Proof of the union counterexample}
\begin{proof}[Proof of Theorem~\ref{union-thm:union-failure}]
Choose
\[
 \frac1p<s<\frac14+\frac1{2p}<\frac12
\]
and take $F_s$ from Lemma~\ref{union-lem:gaussian-jensen}. Since $sp>1$,
\[
 \E\sum_{n\ge0}\frac{|\xi_n|^p}{(n+1)^{sp}}
 =\E|\xi_0|^p\sum_{n\ge0}(n+1)^{-sp}<\infty.
\]
Hence $F_s\in\lp$ almost surely. Fix a realization satisfying this and all conclusions of the lemma. Then $W=Z(F_s)$ and $-W=Z(F_s(-z))$ are disjoint exact zero sets for $\lp$, with simple zeros.

Suppose a nonzero $h\in\lp$ vanishes on their union. Dividing $h$ by its first nonzero monomial preserves membership in $\lp$ and the prescribed zeros, so we may assume $h(0)\ne0$. The two zero sets are disjoint, and all zeros of $F_s$ and $F_s(-\cdot)$
are simple. Since reflection preserves zero moduli, the Jensen sum of $h$
contains the contributions from both zero sets, with no overlap. Hence
\[
 J_h(r_k)\ge J_{F_s}(r_k)+J_{F_s(-\cdot)}(r_k)
 =2J_{F_s}(r_k).
\]
As $r_k=1-e^{-k^2}$, Lemma~\ref{union-lem:jensen-upper} gives
\[
 \limsup_{k\to\infty}\frac{J_h(r_k)}{k^2}\le\frac12-\frac1p,
\]
whereas the almost-sure limit in Lemma~\ref{union-lem:gaussian-jensen}
and the lower bound above imply
\[
 \liminf_{k\to\infty}\frac{J_h(r_k)}{k^2}
 \ge 2\left(\frac12-s\right)>\frac12-\frac1p.
\]
This is a contradiction. The same comparison applies to any two rotations
of this fixed $F_s$ whose zero sets are disjoint, since rotations preserve
all zero moduli.
\end{proof}

\begin{remark}[A Gaussian construction for every $p>2$]\label{rem:gaussian-construction}
Fix $p>2$ and choose $s$ with
\[
 \frac1p<s<\frac14+\frac1{2p}.
\]
For independent standard complex Gaussian random variables $\xi_n$, set
\[
 F_s(z)=\sum_{n\ge0}\frac{\xi_n}{(n+1)^s}z^n,
 \qquad W_+=Z(F_s),\qquad W_-=Z(F_s(-z))=-W_+.
\]
With probability one, $F_s$ and $F_s(-z)$ belong to $\ell_A^p$ and
have only simple zeros; their zero sets $W_+$ and $W_-$ are disjoint
exact zero sets for $\ell_A^p$, and $W_+\cup W_-$ is a uniqueness
set for $\ell_A^p$.
These assertions follow from the construction in the proof of
Theorem~\ref{union-thm:union-failure}: the inequality $sp>1$ gives
coefficient summability, while Lemmas~\ref{union-lem:jensen-upper}
and~\ref{union-lem:gaussian-jensen} give the strict comparison
\[
 2\left(\frac12-s\right)>\frac12-\frac1p.
\]
For instance, when $p=4$, one may take $s=1/3$ and
\[
 F(z)=\sum_{n\ge0}\frac{\xi_n}{(n+1)^{1/3}}z^n.
\]
Almost surely, $F,F(-z)\in\ell_A^4$ have only simple zeros,
$Z(F)\cap Z(F(-z))=\varnothing$, and $Z(F)\cup Z(F(-z))$ is a
uniqueness set for $\ell_A^4$: indeed $1/4<1/3<3/8$, and the Jensen
comparison is $2\cdot(1/6)=1/3>1/4$.
\end{remark}

\section{Application to shift-invariant subspaces}\label{sec:invariants}
The construction in Section~\ref{sec:unions} yields more than a single
pair of zero sets. We first obtain a family of the cardinality of the continuum for
which every pair has the same uniqueness property, and then pass to the
corresponding vanishing subspaces. Throughout this section,
$\mathfrak c=|\mathbb R|$ denotes the cardinality of the continuum.

\begin{proposition}[A continuum of rotated zero sets]
\label{prop:continuum-zero-sets}
For every $2<p<\infty$, there are $F\in\ell_A^p$ with $F(0)\ne0$ and
only simple zeros, and a set $\Lambda\subset\T$ of cardinality
$\mathfrak c$, containing $1$, such that the sets
\[
 W_\lambda=\lambda Z(F),\qquad \lambda\in\Lambda,
\]
are pairwise disjoint exact zero sets for $\ell_A^p$, and
$W_\lambda\cup W_\mu$ is a uniqueness set for $\ell_A^p$ whenever
$\lambda\ne\mu$. Each $W_\lambda$ is realized exactly by
$F_\lambda(z)=F(\overline\lambda z)$, with
$\|F_\lambda\|_p=\|F\|_p$.
\end{proposition}
\begin{proof}
Fix $F=F_s$ as in the proof of
Theorem~\ref{union-thm:union-failure}, with
$1/p<s<1/4+1/(2p)$, and put $W=Z(F)$.
The set $W$ is countably infinite and excludes zero; infinitude follows
from the Jensen growth in Lemma~\ref{union-lem:gaussian-jensen}.

The set of rotations that can identify two points of $W$ is
\[
 D_W=\{w'/w:w,w'\in W,\ |w|=|w'|\}\subset\T.
\]
It is countable. Let $H$ be the subgroup of $\T$ generated by $D_W$,
that is, the set of finite products of members of $D_W$ and their inverses.
It is countable as well. Choose $\Lambda$ to contain exactly one
representative of each coset of $H$ in $\T$, choosing $1$ for the coset
$H$. Since every coset is countable and $|\T|=\mathfrak c$, we have
$|\Lambda|=\mathfrak c$. If $\lambda w=\mu w'$ for
$\lambda,\mu\in\Lambda$ and $w,w'\in W$, then
$\lambda/\mu=w'/w\in D_W\subset H$. The two representatives therefore
coincide. Thus the sets $W_\lambda$ are pairwise disjoint.

Rotation gives $Z(F_\lambda)=W_\lambda$,
$\|F_\lambda\|_p=\|F\|_p$, and $J_{F_\lambda}=J_F$.
For $\lambda\ne\mu$, repeat the final comparison in the proof of
Theorem~\ref{union-thm:union-failure}, with $F_\lambda,F_\mu$ in place
of $F_s,F_s(-\cdot)$. Dividing out a power of $z$ if necessary, assume $h(0)\ne0$.
The two zero sets are disjoint, and $F_\lambda$ and $F_\mu$ have only
simple zeros. Since rotations preserve zero moduli, the Jensen sum of $h$
contains the contributions from both zero sets, with no overlap. Thus
\[
 J_h(r)\ge J_{F_\lambda}(r)+J_{F_\mu}(r)=2J_F(r).
\]
But the choice $s<\frac14+\frac1{2p}$ gives
$2(\frac12-s)>\beta_p$, contradicting the upper bound in
Lemma~\ref{union-lem:jensen-upper}.
The argument is deterministic after $F$ has been fixed and applies to
every pair; no uncountable intersection of probability-one events is used.
\end{proof}

Let $S$ be the right shift on the complex sequence space
$\ell^p(\mathbb N_0)$:
\[
 S(a_0,a_1,a_2,\ldots)=(0,a_0,a_1,\ldots).
\]
A linear subspace $M$ is \emph{shift-invariant} if $SM\subset M$.
Under the isometric identification
$(a_n)_{n\ge0}\longleftrightarrow f(z)=\sum_{n\ge0}a_nz^n$,
$S$ is multiplication by $z$ on $\ell_A^p$.
The existence, for $p>2$, of a family of cardinality continuum of nonzero
closed shift-invariant subspaces with pairwise trivial intersections was
established in~\cite{Ab}.
The following corollary gives a realization by zero conditions.

\begin{corollary}[A continuum of invariant subspaces with pairwise trivial intersections]
\label{cor:disjointshift}
For every $2<p<\infty$, there is a family
$\{M_\lambda:\lambda\in\Lambda\}$, $|\Lambda|=\mathfrak c$, of
nonzero closed, infinite-dimensional subspaces of
$\ell^p(\mathbb N_0)$ such that
\[
 SM_\lambda\subset M_\lambda,\qquad
 M_\lambda\cap M_\mu=\{0\}\quad(\lambda\ne\mu).
\]
The subspaces may be chosen as diagonal rotations of a single subspace:
with $1\in\Lambda\subset\T$,
\[
 M_\lambda=R_\lambda M_1,\qquad
 (R_\lambda a)_n=\overline\lambda^{\,n}a_n.
\]
\end{corollary}
\begin{proof}
Work first in $\ell_A^p$ and take $F_\lambda$ and $W_\lambda$ from
Proposition~\ref{prop:continuum-zero-sets}. For $E\subset\D$, write
\[
 I_p(E)=\{f\in\ell_A^p:f(w)=0\text{ for every }w\in E\}.
\]
This \emph{vanishing subspace} is closed because each point evaluation
is continuous by~\eqref{eq:evaluation}. It is invariant under
multiplication by $z$, since $(zf)(w)=wf(w)$.

Set $M_\lambda=I_p(W_\lambda)$. Each subspace is nonzero, since it
contains $F_\lambda$. For distinct parameters, the uniqueness property
of their union gives
\[
 M_\lambda\cap M_\mu
 =I_p(W_\lambda\cup W_\mu)=\{0\}.
\]
The functions $F_\lambda,zF_\lambda,z^2F_\lambda,\ldots$ belong to
$M_\lambda$ and are linearly independent: a finite linear relation
would give $P F_\lambda\equiv0$ for a polynomial $P$, hence $P=0$.
Thus every $M_\lambda$ is infinite-dimensional. In particular the
subspaces are distinct, so the family has cardinality $\mathfrak c$.
Finally, $f(z)\mapsto f(\overline\lambda z)$ maps $I_p(W_1)$ onto
$I_p(W_\lambda)$ isometrically. Passing to coefficient sequences gives
the stated diagonal rotations and the result in $\ell^p(\mathbb N_0)$.
\end{proof}

The same conclusion holds for the a priori smaller cyclic subspaces
\[
 N_\lambda=\overline{\operatorname{span}}^{\,\ell_A^p}
   \{z^nF_\lambda:n\ge0\}\subset M_\lambda,
\]
which are nonzero and infinite-dimensional and inherit pairwise trivial
intersections from the $M_\lambda$.

\end{document}